\documentclass[10pt,twocolumn,twoside,journal]{IEEEtran} 

\usepackage[dvipsnames]{xcolor}
\usepackage{cite}
\usepackage{amsmath} 
\usepackage{amssymb}  
\usepackage{lipsum}
\usepackage{amsfonts}
\usepackage{dblfloatfix}% http://ctan.org/pkg/dblfloatfix
\usepackage{mathrsfs}
\usepackage{mathtools}
\usepackage[dvipsnames]{xcolor}
\usepackage[thmmarks, amsmath]{ntheorem}
\usepackage[shortlabels]{enumitem}
\usepackage[bbgreekl]{mathbbol}
\usepackage{multirow}
\usepackage{booktabs}
\usepackage{pdfcomment}
\hypersetup{hidelinks}

\DeclareSymbolFontAlphabet{\mathbb}{AMSb}
\DeclareSymbolFontAlphabet{\mathbbl}{bbold}

\DeclareMathOperator{\diag}{diag}

\DeclareMathOperator{\diff}{d}

\providecommand{\norm}[1]{\lVert#1\rVert}
\providecommand{\norm}[1]{\lvert#1\rvert}
\newcommand{\R}{{\mathbb R}}
\newcommand{\N}{{\mathbb N}}
\newcommand{\mc}{\mathcal}
\newcommand{\ddt}{\tfrac{\diff}{\diff \!t}}
\newcommand{\dxdy}[2]{\frac{\partial #1}{\partial #2}}

\newcommand{\Yn}{Y_{\text{\upshape{net}}}}

\newtheorem{theorem}{Theorem}
\newtheorem{lemma}{Lemma}
\newtheorem{proposition}{Proposition}
\newtheorem{assumption}{Assumption}
\newtheorem{remark}{Remark}
\newtheorem{definition}{Definition}
\newtheorem{condition}{Condition}
\newtheorem{corollary}{Corollary}

\allowdisplaybreaks

\begin{document}

\title{\LARGE \bf  A Lyapunov framework for nested dynamical systems on multiple time scales with application to converter-based power systems \thanks{This work was partially funded by the Swiss Federal Office of Energy under grant number SI/501707.}}
\author{Irina Suboti\'c, Dominic Gro\ss{}, Marcello Colombino, and  Florian D\"orfler  \thanks{I. Suboti\'c, M. Colombino, and F. D\"orfler, are with the Automatic Control Laboratory at ETH Z\"urich, Switzerland, D. Gro\ss{} is with the Department of Electrical and Computer Engineering at the University of Wisconsin-Madison, USA; e-mail:\{subotici,mcolombi,dorfler\}@ethz.ch,  dominic.gross@wisc.edu}}

\maketitle
\begin{abstract}
In this work, we present a Lyapunov framework for establishing stability with respect to a compact set for a nested interconnection of nonlinear dynamical systems ordered from slow to fast according to their convergence rates, where each of the dynamics are influenced only by the slower dynamics and the successive fastest one. The proposed approach explicitly considers more than two time scales, it does not require modeling multiple time scales via scalar time constants, and provides analytic bounds that make ad-hoc time-scale separation arguments rigorous. Motivated by the technical results, we develop a novel control strategy for a grid-forming power converter that consists of an inner cascaded two-degree of freedom controller and dispatchable virtual oscillator control as a reference model. The resulting closed-loop converter-based AC power system is in the form of a nested system with multiple time scales. We apply our technical results to  obtain explicit bounds on the controller set-points, branch powers, and control gains that guarantee almost global asymptotic stability of the multi-converter AC power system with respect to a pre-specified solution of the AC power-flow equations. Finally, we validate the performance of the proposed control structure in a case study using a high-fidelity  simulation with detailed hardware validated converter models.
\end{abstract}

\section{Introduction}

\IEEEPARstart{T}{ime-scale} separation arguments are ubiquitous in control design and analysis of large-scale engineering systems that contain dynamics on multiple time scales from different physical domains. 
Traditionally, singular perturbation theory has been the standard tool to analyze nonlinear dynamics that evolve on multiple time scales \cite{KOS76,K15,PKC82}. Within this framework, stability conditions are typically provided for hyperbolic fixed points of systems with two time scales and a ``small'' scalar time constant describing the fast time scale. The results can be extended to linear systems with two time scales and a fast time scale modeled by multiple time constants \cite{KK79}, slow-fast control systems with non-hyperbolic fixed points \cite{JSP19}, and multiple time scales by successively grouping them into two time scales (see, e.g., \cite{K15}). In contrast, our approach explicitly considers multiple time scales, stability with respect to a compact set, does not require modeling time scales via scalar time constants, and exploits the nested structure typically exhibited by systems with multiple time scales such as power systems \cite{PWS-MAP:98} and biological system \cite{MR92}. To this end, we develop a general Lyapunov function framework for stability analysis of nested nonlinear dynamical systems that can be ordered from slow to fast in terms of convergence rates to their set of steady-states and only depend on the states of slower systems and the next fastest one.

The analysis in this paper is based on a recently developed Lyapunov characterization of almost global asymptotic stability with respect to a compact set presented in \cite{line_dynamics} that requires that the set of states that are unstable but attractive has zero Lebesgue measure. 
The technical contribution is twofold: first, we develop a Lyapunov characterization of unstable hyperbolic fixed points that have a region of attraction of measure zero; then we provide a Lyapunov framework that results in conditions under which the guarantees obtained by applying the aforementioned Lyapunov conditions in \cite{line_dynamics} to a reduced-order system translate to the full-order nested dynamical system. Our results can be interpreted as an extension of the conditions for systems with two time scales in \cite{SK84},\cite[Ch 11]{K02} to multiple nested systems and a more general notion of stability. Moreover, we reduce conservatism by allowing for a wider range of comparison functions. 

Motivated by the transition of power systems towards renewable energy sources that are connected to the system via power electronics \cite{FM-FD-GH-DH-GV:18, WEB+15}, we apply our technical results to multi-converter AC power systems. The analysis and control of power systems and microgrids is typically based on reduced-order models of various degrees of fidelity that exploit the pronounced time-scale separation between the dynamics of synchronous machines and power converters and the transmission network \cite{PWS-MAP:98, PKC82, ADS18, SZO+16, LD14}. 
While these ad-hoc model simplifications have proved themselves useful their validity for converter-based systems is questionable. For instance, the assumption that the dynamics of transmission lines can be neglected breaks down for power systems dominated by fast acting power converters \cite{line_dynamics,vorobev2017high,MSV+19,Guerrero2015}. In this work, we make the time-scale separation argument rigorous by explicitly considering the interaction of dynamics on different time scales (i.e., converter dynamics, inner controls, line dynamics) and quantifying the parameters (e.g., set-points, control gains, transmission line parameters, etc.) for which stability for the overall system can be ensured.

The prevalent approach to grid-forming control is so-called droop-control \cite{MCC-DMD-RA:93,Guerrero2015} and synchronous machine emulation~\cite{zhong2011synchronverters,JAD18,DSF15}. However, while droop control and machine emulation can provide useful insights, stability guarantees for multi-converter systems are typically local and don't extend to line dynamics, detailed converter models, and operating points with non-zero power flows. Almost global stability guarantees are available for single machines connected to an infinite bus \cite{SEO+2016} and droop controlled multi-converter systems when neglecting the network dynamics and losses, internal converter dynamics and controls, and voltage magnitude dynamics \cite{SEO2019}.

In contrast, virtual oscillator control (VOC) ensures almost global synchronization \cite{johnson2014synchronization,MS-FD-BJ-SD:14b,LABT-JPH-JM:13} but cannot be dispatched, and almost global asymptotic stability with line dynamics can be ensured for dispatchable virtual oscillator control (dVOC) \cite{colombino2017global2,GC+19} for appropriate control gains and power converters modeled as controllable voltage sources \cite{line_dynamics}. All of these works use simplified models that can be regard as the slow time-scales of a nested dynamical system. One contribution of this work is to illustrate how to include faster time scales.

We consider a converter-based power system model that includes dynamic models of voltage source converters as well as transmission network dynamics. We develop a control strategy for the converters that uses dispatchable virtual oscillator control (dVOC) as a reference model for a cascaded two-degree of freedom voltage and current controller in stationary $\alpha\beta$ coordinates instead of local $dq$ coordinates. This allows us to model the whole inverter based power system as a nested interconnection of subsystems (dVOC, the transmission line dynamics, the inner control loops) that evolve on different time-scales and to apply our novel Lyapunov-based stability criterion to provide explicit stability conditions and show that there always exists parameters (i.e., control gains and set-points) that guarantee almost global asymptotic stability of the overall system. Finally, we validate the proposed control architecture in a high-fidelity simulation of the hardware set-up described in \cite{GC+19}.

The remainder of this section recalls some basic notation and results from graph theory. Section \ref{sec.lyap} provides definitions and a preliminary technical results on almost global asymptotic stability with respect to sets. The main theoretical contribution is given in Section \ref{section:SPT.extension}. In Section \ref{sec:model} we present a detailed model of a multi-converter AC power system and the control objectives. Section \ref{sec:control_design} presents a cascaded two-degree of freedom control structure that tracks a reference obtained by dVOC as a reference model and Section \ref{section:Stability.Analysis} presents stability conditions for the multi-converter system. The results are illustrated using a high-fidelity simulation study in Section~\ref{sec:illustrative.example}, and Section~\ref{sec.conclusion} provides the conclusions.

\subsection*{Notation}
We use $\R$ and $\mathbb N$ to denote the set of real and natural numbers and define $\R_{\geq a}\coloneqq \{x \in \mathbb R \vert x \geq a\}$ and, e.g., $\R_{[a,b)}\coloneqq \{x \in \mathbb R \vert a \leq x < b\}$. Given ${\theta\in[-\pi,\pi]}$ the 2D rotation matrix is given by
\begin{align*}
\mc R(\theta) \coloneqq \begin{bmatrix}\cos(\theta) & -\sin(\theta)\\  \sin(\theta) & \cos(\theta) \end{bmatrix} \in \R^{2 \times 2}.
\end{align*}
Moreover, we define the $90^\circ$ rotation matrix $J \coloneqq \mc R(\pi/2)$ that can be interpreted as an embedding of the complex imaginary unit $\sqrt{-1}$ into $\R^2$. Given a matrix $A$, $A^\mathsf{T}$ denotes its transpose. We use $\|A\|$ to indicate the induced 2-norm of $A$. We write $A\succcurlyeq0$  $(A\succ0)$ to denote that $A$ is symmetric and positive semidefinite (definite). For column vectors $x\in\R^n$ and $y\in\R^m$ we use $(x,y) = [x^\mathsf{T}, y^\mathsf{T}]^\mathsf{T} \in \R^{n+m}$ to denote a stacked vector, and $\norm{x}$ denotes the Euclidean norm. The absolute value of a scalar $y \in \R$ is denoted by $|y|$. Furthermore, $I_n$ denotes the identity matrix of dimension $n$, and $\otimes$ denotes the Kronecker product. For any matrix $M$ and any $n\in \N$, we define $M_n = I_n \otimes M$. Matrices of zeros of dimension $n \times m$ are denoted by $\mathbbl{0}_{n\times m}$ and $\mathbbl{0}_{n}$ denotes column vector of zeros of length $n$. We use $\norm{x}_C \coloneqq \min_{z \in \mathcal{C}} \norm{z-x}$ to denote the distance of a point $x$ to a set $\mathcal{C}$. We use $\varphi_f(t,x_0)$ to denote the solution of $\ddt x = f(x)$ at time $t \geq 0$ starting from the initial condition $x(0)=x_0$ at time $t_0=0$.

\section{Almost global asymptotic stability with respect to a compact set}\label{sec.lyap}

Consider the dynamical system
\begin{align}\label{eq:dynsys}
 \ddt x = f(x),
\end{align}
where $x \in \mathbb{R}^n$ denotes the state vector and $f: \mathbb{R}^n \to \mathbb{R}^n$ is a Lipschitz continuous function. In order to state the main results of the paper, we require the following definition of almost global asymptotic stability with respect to a set \cite{A04}.
\begin{definition}{\bf(Almost global asymptotic stability)}\label{def:ags}
	The dynamic system \eqref{eq:dynsys} is called almost globally asymptotically stable with respect to a compact set $\mc C\subset \mathbb{R}^n$ if 
	\begin{enumerate}[label=(\roman*)]
		\item it is almost globally attractive with respect to $\mc C$, i.e.,
		\begin{align}
		\lim_{t\to\infty} \norm{\varphi_f(t,x_0)}_{\mc C} = 0 \label{eq.AGA}
		\end{align}
		holds for all $x_0 \notin \mc Z$, and $\mc Z$ has zero Lebesgue measure, 
		\item it is Lyapunov stable with respect to $\mc C$, i.e., for every $\varepsilon \in \mathbb{R}_{>0}$ there exists $\delta \in \mathbb{R}_{>0}$ such that
		\begin{align}\label{eq:lyapstab}
		\norm{x_0}_{\mc C} < \delta \implies \norm{\varphi_f(t,x_0)}_{\mc C} < \varepsilon, \qquad \forall t \geq 0.
		\end{align}
	\end{enumerate}
\end{definition}
Next, we recall the definition of comparison functions used to establish stability properties of dynamical systems \cite{H67}.
\begin{definition}{\bf(Comparison functions)}
	A function $\chi_c:\mathbb{R}_{\geq0} \to \mathbb{R}_{\geq0}$ is of class $\mathscr{K}$ if it is continuous, 
	strictly increasing and $\chi_c(0)=0$; it is
	of class  $\mathscr{K}_{\infty}$ if it is a $\mathscr{K}$-function and $\chi_c(s)\to\infty$ as $s\to\infty$. 
\end{definition}
Next, consider a set $\mc U \subset \mathbb{R}^n$ that is invariant with respect to \eqref{eq:dynsys} (i.e., $\varphi_f(t,x_0) \in \mc U$ for all $t \in \R_{>0}$ and all $x_0 \in \mc U$), satisfies $\mc C \cap \mc U = \emptyset$, and 
corresponds to e.g., undesirable equilibria or limit cycles of \eqref{eq:dynsys}. In this case global asymptotic stability of \eqref{eq:dynsys} with respect to $\mc C$ cannot be established. Instead, the following Theorem provides a Lyapunov function characterization of almost global asymptotic stability with respect to $\mc C$ \cite[Th. 1]{line_dynamics}.
\begin{theorem} {\bf(Lyapunov functions)}\label{theorem:AGAS}
	Consider a compact set $\mc C \subset \mathbb{R}^n$ and a zero Lebesgue measure set $\mc U \subset \mathbb{R}^n$ that is invariant with respect to \eqref{eq:dynsys}. Moreover, consider a continuously differentiable function $\mc V: \mathbb{R}^n \to \mathbb{R}_{>0}$ and comparison functions $\chi_1, \chi_2 \in \mathscr{K}_\infty$ and $\chi_3 \in \mathscr{K}$ such that 
	\begin{align*}
	\chi_1(\norm{x}_{\mathcal{C}}) \leq \mc V(x) &\leq \chi_2(\norm{x}_{\mathcal{C}})\\
	\frac{\diff}{\diff t}  \mc V(x) \coloneqq \frac{\partial \mc V}{\partial x} f(x) &\leq -\chi_3(\norm{x}_{\mc C \cup \, {\mc U}})
	\end{align*} 
	holds for all $x \in \mathbb{R}^n$. Moreover, let 
	\begin{align*}
	\mc Z_{\mc U,f} \coloneqq \{ x_0 \in \mathbb{R}^n \vert \lim\nolimits_{t\to\infty} \norm{\varphi_{f}(t,x_0)}_{\mc U} = 0 \}
	\end{align*}
	denote the region of attraction of $\mc U$ under \eqref{eq:dynsys}. If  $\mc Z_{\mc U,f}$ has zero Lebesgue measure, the dynamics \eqref{eq:dynsys} are almost globally asymptotically stable with respect to $\mathcal{C}$. 
\end{theorem}
Conceptually Theorem \ref{theorem:AGAS} is similar to the main result in \cite{AE2015} that has been used in \cite{SEO+2016} to study almost global stability of a synchronous machine connected to an infinite bus. However, \cite[Def. 8]{AE2015} uses the bound $\chi_1(\norm{x}_{\mc W}) \leq \mc V(x) \leq \chi_2(\norm{x}_{\mc W}+c)$ with $\mc W = \mc C \cup \mc U$ and $c \in \R_{\geq 0}$ and establishes practical stability of $\mc W$. In contrast, the Lyapunov stable set $\mc C$ and the globally attractive set $\mc C \cup \, \mc U$ are treated separately in Theorem \ref{theorem:AGAS} to obtain a more precise characterization of stability.

Besides finding a suitable Lyapunov function, the main difficulty in applying Theorem \ref{theorem:AGAS} is to verify that the region of attraction $\mc Z_{\mc U,f}$ of the (undesirable) attractive set $\mc U$, has measure zero. To this end, our first contribution is a Lyapunov-like condition that characterizes unstable hyperbolic fixed points with region of attraction that has zero Lebesgue measure.
\begin{theorem} \label{theorem:ROA} {\bf{(Region of attraction)}}		Consider a continuously differentiable function $f: \mathbb{R}^n \to \mathbb{R}^n$ and a fixed point  $x^\star \in \mathbb{R}^n$ (i.e., $f(x^\star)=\mathbbl{0}_N$). Moreover, consider the Jacobian $A_{x^\star} \coloneqq \frac{\partial f(x)}{\partial x}\big\vert_{x=x^\star}$, the linearized dynamics $\ddt \delta_x = A_{x^\star} \delta_x$, and a quadratic function $\mc V_\delta \coloneqq \frac{1}{2}\delta^\mathsf{T}_x P_\delta \delta_x$ with $P_\delta = P^\mathsf{T}_\delta \in \mathbb{R}^n$. If $\ddt \mc V_\delta (\delta_x) \coloneqq \frac{ \partial \mc V_\delta }{\partial x_\delta} A_{x^\star} \delta_x < 0$ holds for all $\delta_x \neq \mathbbl{0}_n$ and there exists $\delta_x^\prime \in \mathbb{R}^n$ such that $ \mc V_\delta (\delta^\prime_x) <  \mc V_\delta (\mathbbl{0}_n)$, then the region of attraction $\mc Z_{\{x^\star\},f}$ of $x^\star$ under the nonlinear dynamics $\ddt x = f(x)$ has Lebesgue measure zero.
\end{theorem}
\begin{IEEEproof}
By \cite[Prop. 4.1]{PM09} the region of attraction $\mc Z_{\{x^\star\},f}$ of $x^\star$ under the nonlinear dynamics $\ddt x = f(x)$ has Lebesgue measure zero if $f: \mathbb{R}^n \to \mathbb{R}^n$ is continuously differentiable and the Jacobian $A_{x^\star} $ has at least one eigenvalue with positive real part. To show that at least one eigenvalue of $A_{x^\star}$ has positive real part, we first note that $\mc V_\delta (\delta^\prime_x) < \mc V_\delta (\mathbbl{0}_n)$ and $\ddt \mc V_\delta   < 0$ imply that $\mc V_\delta (\varphi_{A_{x^\star}}(t,c^\prime \delta^\prime_x)) \to -\infty$ and $\norm{\varphi_{A_{x^\star}}(t,c^\prime \delta^\prime_x)} \to \infty$ as $t \to \infty$ for any $c^\prime \in \R_{>0}$, i.e., $A_{x^\star}$ is not stable. Next, assume that $A_{x^\star}$ has eigenvalues with zero real part. Then there exist initial conditions $\delta_{x_0} \neq \mathbbl{0}_n$ such that $\varphi_{A_{x^\star}}(t,\delta_{x_0})$ remains bounded for all $t \in \mathbb{R}_{\geq 0}$, but does not converge to the origin $\mathbbl{0}_n$. However, because $\mc V_\delta (\varphi_{A_{x^\star}}(t,\delta_{x_0}))$ is strictly decreasing in $t$ for all $\delta_{x_0} \neq \mathbbl{0}_n$ it either holds that $\norm{\varphi_{A_{x^\star}}(t,\delta_{x_0})} \to 0$ as $t \to \infty$ or $\mc V_\delta (\varphi_{A_{x^\star}}(t,\delta_{x_0})) \to -\infty$ and $\norm{\varphi_{A_{x^\star}}(t,\delta_{x_0})} \to \infty$ as $t \to \infty$. Therefore, $A_{x^\star}$ cannot have eigenvalues with zero real part and, because $A_{x^\star}$ is unstable, at least one eigenvalue must have positive real part and $\mc Z_{\{x^\star\},f}$ has Lebesgue measure zero.
\end{IEEEproof}
Note that the conditions and steps that establish instability of the fixed-point $x^\star$ are a special case of Chetaev's theorem \cite[Th. 4.3]{K02}. However, Chetaev's theorem does not establish that at least one eigenvalue of $A_{x^\star}$ has positive real part.

\section{Stability Theory for Nested Systems on multiple time scales} \label{section:SPT.extension}
In this section we will present a model of nested dynamical systems on multiple time-scales and extend the results from Theorem \ref{theorem:AGAS} and Theorem \ref{theorem:ROA} to this class of systems.

\subsection{Nested systems on multiple time scales}
\begin{figure}[t!!!]
\centering
\includegraphics[width=.98\columnwidth]{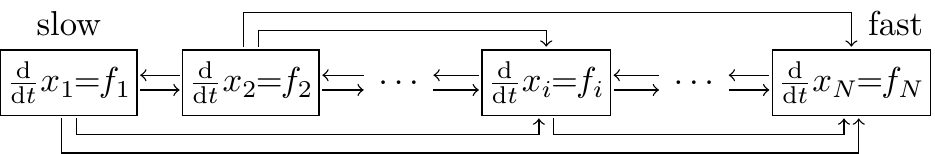}
\vspace{.2em}
\caption{Block diagram of the nested dynamical system on $N$ time scales \label{fig:block.diagram}}
\vspace{2em}
\includegraphics[width=.98\columnwidth]{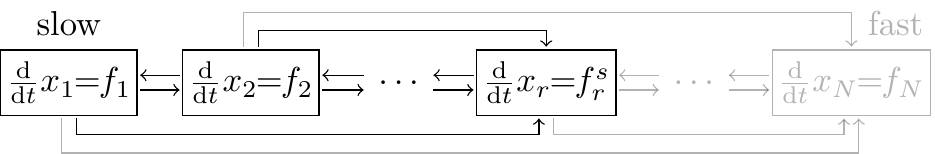}
\vspace{.2em}
\caption{Block diagram of the reduced-order dynamics on $r$ time scales (black) and neglected dynamics (grey). \label{fig:block.diagram.red}}
\end{figure}

Consider the nested dynamical system on $N$ time scales shown in Figure \ref{fig:block.diagram} given by
\begin{subequations} \label{eq.General.N.dimension}
	\begin{align}
		\ddt x_i &=  f_i (x_1, \hdots, x_{i+1}),  \ \forall i \in \N_{[1,N-1]}, \\
		\ddt x_N &= f_N (x_1, \hdots, x_{N}), 
	\end{align}
\end{subequations}
where $x_i \in \R^{n_i}$ denotes the state of the subsystem $i \in \N_{[1,N-1]}$ and $f_i:\R^{n_1} \times \hdots \times \R^{n_{i+1}} \to \R^{n_i}$ and $f_N:\R^{n_1} \times \hdots \times \R^{n_{N}} \to \R^{n_N}$ are Lipschitz continuous functions. Letting $x=(x_1,\ldots,x_N)$ and $f=(f_1,\ldots,f_N)$, we obtain the full system dynamics \eqref{eq:dynsys} with $n=\sum_{i=1}^N n_i$.
Broadly speaking, we assume that the dynamics are ordered from slow to fast convergence to their set of steady-states, i.e., the outer dynamics are the slowest and the inner right dynamics are the fastest (see Figure \ref{fig:block.diagram}). To make this argument precise, we recursively define the steady-state maps for each $f_i$, $i \in \mathbb{N}_{[2,N]}$.
{
\begin{assumption}{\bf(Steady-state maps)}\label{assumption:steadystate}\\
We assume that there exists a unique steady-state map $\phi_N(x_1,\ldots,x_{N-1})$ such that, for all $x_1,...,x_{N-1}$, we have that $f_N(x_1,\ldots,x_{N-1},\phi_N(x_1,\ldots,x_{N-1}))=0$. In a similar fashion,  we recursively define for $i\in \N_{[1,N-1]}$ the maps
\begin{multline*}
f^s_{N-i}(x_1,\ldots,x_{N-i}) \coloneqq \\ f_{N-i}(x_1,\ldots,x_{N-i},\phi_{N-i+1}(x_1,\ldots,x_{N-i})),
\end{multline*}
and recursively assume that, for all $i \in \N_{[1,N-2]}$, there exists a unique steady-state map $\phi_{N-i}(x_1,\ldots,x_{N-i-1})$ such that $f^s_{N-i}(x_1,\ldots,x_{N-i-1},\phi_{N-i}(x_1,\ldots,x_{N-i-1}))=0$ for all $x_1,\ldots,x_{N-i-1}$.
\end{assumption}}
{We emphasize that $\phi_i(x_1,...,x_{i-1})$ is a function of time-varying argument and only describes a constant steady-state of the $i^{\text {th}}$ subsystem if the ``slower'' states $x_1,...,x_{i-1}$ are constant. Moreover,} $f^s_i$ denotes the vector field corresponding to the dynamics with index $i \in \N_{[1,N-1]}$ with the state $x_{i+1}$ restricted to its steady-state map, i.e., $x_{i+1}=\phi_{i+1}(x_1,...,x_i)$. Intuitively, if the systems are ordered from slow to fast convergence to their set of steady-states, Assumption \ref{assumption:steadystate} suggests the natural model reduction procedure that successively replaces fast dynamics by their steady-state maps. Given $r \in \N_{[1,N-1]}$, this results in the reduced-order dynamics  (see Figure \ref{fig:block.diagram.red}) with state vector $(\tilde{x}_1,\ldots,\tilde{x}_r) \in \R^{\sum_{i=1}^r n_i}$ given by
\begin{subequations}\label{eq:redordermodel}
\begin{align}
    \ddt \tilde{x}_i &=  f_i (\tilde{x}_1, \hdots, \tilde{x}_{i+1}),  \ &&\forall i \in \N_{[1,r-1]}, \label{eq:redordermodel:a}\\ 
	\ddt \tilde{x}_r &= f^s_r (\tilde{x}_1, \hdots, \tilde{x}_r),\\
	\tilde{x}_i &= \phi_i(\tilde{x}_1,\ldots,\tilde{x}_{i-1}), \ &&\forall i \in \N_{[r+1,N]}.
\end{align}
\end{subequations}
We emphasize that \eqref{eq:redordermodel} defines a reduced-order dynamical system. {In particular, for $r=1$ all fast dynamics are neglected, \eqref{eq:redordermodel:a} is not defined, and we obtain the reduced-order dynamics $\ddt \tilde{x}_1 = f^s_1(\tilde{x}_1)$ of order $n_1$.} Next, given a set $\Omega \subset \R^{n_1}$, we define the mapping of $\Omega$ under the steady-state maps as
\begin{align}\label{eq.Xs}
\!\!\mc X^s(\Omega) \coloneqq \left\{ \begin{bmatrix} x_1 \\ \vdots \\ x_N \end{bmatrix} \in \R^{n} \left\vert\; \!\!\! \begin{array}{c}
x_1 \in \Omega \\ x_2=\phi_2(x_1) \\  \vdots \\
x_N = \phi_N(x_1,\ldots,x_{N-1})
\end{array}  \right. \!\! \right\}.
\end{align}
In the remainder of this section we derive conditions that allow us to extend guarantees of the type given in Theorem \ref{theorem:AGAS} and Theorem \ref{theorem:ROA} for the reduced-order system $\ddt \tilde{x}_1 = f^s_1(\tilde{x}_1)$ and a set $\mc C_1 \subset \R^{n_1}$ to the full dynamics \eqref{eq.General.N.dimension} and $\mc X^s(\mc C_1) \subset \R^n$.

For ease of notation, when they are obvious by the context, we will drop the arguments of $\phi_i$ in the remainder of the paper.

\subsection{Lyapunov function for nested systems}
For all $i \in \N_{[2,N]}$ we use $y_i \coloneqq x_i-\phi_i \in \R^{n_i}$ to denote the difference of $x_i$ to its steady-state map $\phi_i$. For all $i \in \N_{[1,N]}$ we define the  continuously differentiable Lyapunov function candidates $V_i: \R^{n_i}\to \R_{\geq 0}$. Given positive constants $\mu_i \in \R_{> 0}$ to be determined, a Lyapunov function candidate $\nu : \R^{n_1} \times \ldots \times \R^{n_N} \to \R_{\geq 0}$ for the system \eqref{eq.General.N.dimension} is given by
\begin{align}\label{eq.muLyap}
 \nu \coloneqq \mu_1 V_1(x_1) + \sum\nolimits_{i=2}^{N} \mu_i V_i (y_i).
\end{align}
For clarity of the exposition, we omit the arguments of $V_i$, $\phi_i$, $f_i$ in the remainder. We require the following assumption that bounds the decrease of the individual Lyapunov functions $V_i$ in \eqref{eq.muLyap} for their associated reduced-order models and bounds their increase due to neglecting slower and faster dynamics.
\begin{assumption} \label{assumption:derivative}
Given compact sets $\mc C_1 \subset \R^{n_1}$ and $\mc U_1 \subset \R^{n_1}$, for all $i \in \N_{[1,N]}$, there exist positive constants $\alpha_i \in \R_{>0}$ and $\alpha^\prime_i \in \R_{>0}$, positive
semidefinite functions $\psi_i: \R^{n_i} \to \R_{\geq 0}$ and $\psi^\prime_i: \R^{n_i} \to \R_{\geq 0}$, and $\mathscr{K}$ functions $\sigma_i$ such that 
$\sigma_1(\norm{x_1}_{\mc C_1 \cup \mc\, U_1}) \leq \psi_1(x_1)+\psi^\prime_1(x_1)$ and 
$\sigma_i(\norm{y_i}) \leq \psi_i(y_i)+\psi^\prime_i(y_i)$ for all $i \in \N_{[2,N]}$. Moreover,
\begin{align*}\label{eq:lyapdecr1}
\dxdy{V_1}{x_1} f^s_1 &\leq -\alpha_1 \psi_1(x_1)^2 -\alpha^\prime_1 {\psi^\prime_1(x_1)}^2,
\end{align*}
holds, and for all $i \in \N_{[2,N-1]}$ it holds that
\begin{align*}\label{eq:lyapdecri}
\dxdy{V_i}{y_i} f^s_i &\leq -\alpha_i \psi_i(y_i)^2 -\alpha^\prime_i \psi^\prime_i(y_i)^2,\\
\dxdy{V_N}{y_N} f_N &\leq -\alpha_N \psi_N(y_N)^2 -\alpha^\prime_N \psi^\prime_N(y_N)^2.
\end{align*}
Additionally, there exist constants $\beta_{i,i+1} \in \R_{> 0}$ such that
\begin{align*}\label{eq:lyaperr1}
\dxdy{V_1}{x_1} \big (f_1-f^s_1\big) &\leq {\beta_{1,2}} \psi_1(x_1) \psi_2(y_2)
\end{align*}
holds (i.e., for $i=1$), and for all $i \in \N_{[2,N-1]}$ it holds that
\begin{align*}\label{eq:lyaperri}
\dxdy{V_i}{y_i} \big (f_i-f^s_i\big) &\leq {\beta_{i,i+1}} \psi_i(y_i) \psi_{i+1}(y_{i+1}).
\end{align*} 
Finally, for all $i \in \N_{[2,N]}$, $k \in \N_{[1,i-1]}$, and $j \in \N_{[1,k+1]}$, there exists $b_{i,j,k} \in \R$ such that
\begin{align*}\label{eq:lyaperrss}
	-\dxdy{V_i}{y_i} \dxdy{\phi_i}{x_k}  f_k
	\leq  \sum_{j=2}^{k+1}  b_{i,j,k}  \psi_i(y_i) \psi_j(y_j)  +   b_{i,1,k}  \psi_i(y_i) \psi_1(x_1),
\end{align*}
$\beta_{i,1}\coloneqq \sum\nolimits_{k=1}^{i-1}b_{i,1,k} >0$ for all $i \in  \N_{[2,N]}$, and $\beta_{i,j} \coloneqq \sum\nolimits_{k=j-1}^{i-1} b_{i,j,k}>0$ for all $(i,j) \in  \N_{[3,N]} \times \N_{[2,i-1]}$.
\end{assumption}
In particular, the first inequality in Assumption \ref{assumption:derivative} bounds the decrease of the Lyapunov function $V_1$ along the trajectories of the reduced-order model $\ddt x_1 = f^s_1(x_1)$. Moreover, for $i \in \N_{[2,N]}$ the second inequality bounds the decrease of the Lyapunov function candidates $V_i$ in the error coordinates $y_i$ under the assumption that all slower states are constant (i.e., $\ddt x_j = 0$ for all $j<i$) and all faster states are in their steady state (i.e., $x_j = \phi_j$ for all $j>i$). The remaining inequalities bound the additional terms in the time derivative of the Lyapunov function $\nu$ along the full dynamics \eqref{eq.General.N.dimension} that arise because the slower states are generally not constant and the faster states are generally not in their steady state. 

Note that a Lyapunov-based stability proof requires that the right hand side of the first two inequalities in Assumption \ref{assumption:derivative} can be bounded by appropriate comparison functions. If $\psi^\prime_i=0$, then $\psi_i$ needs to be lower bounded by a $\mathscr{K}$-function and satisfy the last three inequalities in Assumption \ref{assumption:derivative}. In contrast, we only require that $\psi_i + \psi^\prime_i$ is bounded from below by a suitable $\mathscr{K}$-function $\sigma_i$ (see, e.g., Theorem \ref{theorem:AGAS}). For instance, when states do not directly affect other time scales they can be excluded in $\psi_i$ and instead included in $\psi^\prime_i$ to obtain improved bounds (see Remark \ref{rem:circcurr} in appendix \ref{app:lyapred}). For simplicity of notation we define $\gamma_i \coloneqq b_{i,i,i-1}$ for all $i \in  \N_{[2,N]}$. Note that $\beta_{i,i-1}$, $\beta_{i-1,i}$, and $\gamma_i$ will be used to bound the difference between the reduced-order model \eqref{eq:redordermodel} with $r=i-1$ and the full dynamics $\ddt x_i = f_i$. Moreover, $\alpha_i$ bounds the convergence rate of $y_i$ under the assumption that all slower states are constant all faster states are in their steady state. 

Moreover, we define $\mu_i$ used in the Lyapunov function \eqref{eq.muLyap} as $\mu_i \coloneqq \prod_{j=1}^{i-1} \frac{\beta_{j,j+1}}{\beta_{j+1,j}}$ for all 
$i \in \mathbb{N}_{[2,N]}$, $\mu_1 = 1$, and we define the symmetric matrix $M$ as follows.
\begin{definition} \label{def.recursive.matrix.M}
Starting from $M_1=\alpha_1$ the symmetric matrix $M \in \R^{N \times N}$ is recursively defined for all $i \in \N_{[2,N]}$ by its leading principal minors
	\begin{align*}%\label{def.M}
		M_i= \left[\begin{array}{c|c}
		M_{i-1} & -\boldsymbol{\beta}_i \mu_i \\[0.3em]
		\hline\\[-1em]
		\star
		& (\alpha_i-\gamma_i) \mu_i
		\end{array} \right],\,		
	\end{align*} 
and $\boldsymbol{\beta}_i \coloneqq (\ldots, \tfrac{1}{2}\beta_{i,i-3}, \tfrac{1}{2}\beta_{i,i-2}, \beta_{i,i-1}) {\in \R^{i-1}_{>0}}$.
\end{definition}

In the next section we show that the derivative of the Lyapunov function $\nu$ along the trajectories of \eqref{eq.General.N.dimension} is bounded by $\ddt \nu \leq -(\psi_1,\ldots,\psi_N)^{\mathsf{T}} M (\psi_1,\ldots,\psi_N)$, i.e., $\nu$ is decreasing if $M$ is positive definite. The main result of this section are two theorems that exploit this fact to establish almost global asymptotic stability of the nested system \eqref{eq.General.N.dimension} using the Lyapunov function $\nu$. Subsequently, we will provide tractable conditions for verifying that $M$ is positive definite (see Section \ref{sec:poscond}).

\subsection{Almost global asymptotic stability of nested systems}
We are now ready to state the main result that establishes almost global asymptotic stability of \eqref{eq.General.N.dimension} with respect to a set $\mc X^s(\mc C_1)$, where $\mc X^s(\cdot)$ is defined in \eqref{eq.Xs}.
\begin{theorem}{\bf(Almost global asymptotic stability of nested systems)} \label{theorem:SPAGAS}
Consider compact sets $\mc C_1 \subset \R^{n_1}$ and $\mc U_1 \subset \R^{n_1}$. Assume that, for all $i\in\N_{[1,N]}$, there exists $\chi^{\scriptscriptstyle{V_i}}_1 \in \mathscr{K}_{\infty}$ and $\chi^{\scriptscriptstyle{V_i}}_2 \in \mathscr{K}_{\infty}$ such that $\chi^{\scriptscriptstyle{V_1}}_1(\norm{x_1}_{\mc C_1}) \leq V_1(x_1) \leq \chi^{\scriptscriptstyle{V_1}}_2(\norm{x_1}_{\mc C_1})$ holds and 
$\chi^{\scriptscriptstyle{V_i}}_1(\norm{y_i}) \leq V_i(y_i) \leq \chi^{\scriptscriptstyle{V_i}}_2(\norm{y_i})$ holds for all $i \in\N_{[2,N]}$. Suppose Assumption \ref{assumption:steadystate} and \ref{assumption:derivative} hold, $M$ is positive definite, and the region of attraction $\mc Z_{\mc X^s(\mc U_1),f}$ of $\mc X^s(\mc U_1)$ has measure zero, then the system \eqref{eq.General.N.dimension} is almost globally asymptotically stable with respect to $\mc X^s(\mc C_1)$.
\end{theorem}

\begin{IEEEproof}
{Let us consider the Lyapunov function candidate $\nu$ defined in~\eqref{eq.muLyap}.} Using Lemma \ref{lemma:decrease.epsilon} (given in Appendix~\ref{app:proofs}) with $\epsilon_i=1$ for all $i \in \N_{[1,N]}$, we obtain {$M=H$, $M^\prime=H^\prime$, and}
\begin{align*}
  \ddt \nu \leq - \left[ \begin{smallmatrix} \psi^{\phantom{\prime}}_1 \\ \psi^{\phantom{\prime}}_2 \\ \vdots \\ \psi^{\phantom{\prime}}_N \end{smallmatrix} \right]^{\mathsf{T}} \!\!{M} \left[ \begin{smallmatrix} \psi^{\phantom{\prime}}_1 \\ \psi^{\phantom{\prime}}_2\\ \vdots \\ \psi^{\phantom{\prime}}_N  \end{smallmatrix} \right] - \left[ \begin{smallmatrix} \psi^\prime_1 \\ \psi^\prime_2 \\ \vdots \\ \psi^\prime_N \end{smallmatrix} \right]^{\mathsf{T}} \!\!{M}^\prime \left[ \begin{smallmatrix} \psi^\prime_1 \\ \psi^\prime_2\\ \vdots \\ \psi^\prime_N  \end{smallmatrix} \right].
\end{align*}
$M^\prime$ is a positive definite diagonal matrix, and, if $M$ positive definite, there exists a positive constant $\alpha_M \in \R_{>0}$ such that $M \succ \alpha_M I_N$ and $M^\prime \succ \alpha_M I_N$. It follows that
\begin{align*}
  \ddt \nu \leq - \alpha_M \sum\nolimits_{i=1}^{N} \psi^2_i + {\psi^\prime_i}^2 \leq - \tfrac{1}{2}\alpha_M \sum\nolimits_{i=1}^{N} (\psi_i + {\psi^\prime_i})^2,
\end{align*}
where the last inequality follows from the fact that  
$- \alpha_M (\psi^2_i + {\psi^\prime_i}^2) \leq - \tfrac{1}{2}\alpha_M (\psi_i + {\psi^\prime_i})^2$ is equivalent to $\tfrac{1}{2}\alpha_M(\psi_i-\psi^\prime_i)^2 \geq 0$ for all $i \in \N_{[1,N]}$. Using Assumption \ref{assumption:derivative}, we conclude that
\begin{align}\label{eq.lyap.comp.decr}
\ddt \nu \leq  - \tfrac{1}{2}\alpha_M \bigg(\sigma_1(\norm{x_1}_{\mc C_1 \cup \mc U_1})^2  +  \sum\nolimits_{i=2}^{N} \sigma_i(\norm{y_i})^2\bigg).
\end{align}
The right hand side of \eqref{eq.lyap.comp.decr} is negative definite w.r.t $\mc X^s(\mc C_1) \cup \mc X^s(\mc U_1)$, and $\mc X^s(\mc C_1) \cup \mc X^s(\mc U_1)$ is compact. Using the same steps as in \cite[p. 98]{H67} there exists a function $\chi_3 \in \mathscr{K}$ such that
\begin{align*}
\ddt \nu \leq -\chi_3(\norm{x}_{\mc X^s(\mc C_1) \cup \mc X^s(\mc U_1)}).
\end{align*}
Moreover, under the hypothesis of the theorem, for all $i\in\N_{[1,N]}$, there exists $\chi^{\scriptscriptstyle{V_i}}_1 \in \mathscr{K}_{\infty}$ and $\chi^{\scriptscriptstyle{V_i}}_2 \in \mathscr{K}_{\infty}$ such that $\chi^{\scriptscriptstyle{V_1}}_1(\norm{x_1}_{\mc C_1}) \leq V_1(x_1) \leq \chi^{\scriptscriptstyle{V_1}}_2(\norm{x_1}_{\mc C_1})$ holds and 
$\chi^{\scriptscriptstyle{V_i}}_1(\norm{y_i}) \leq V_i(y_i) \leq \chi^{\scriptscriptstyle{V_i}}_2(\norm{y_i})$ holds for all $i \in\N_{[2,N]}$. 
For $j \in \N_{[1,2]}$, we define the functions $\tilde{\chi}_j \coloneqq \sum_{i=1}^{N} \mu_i \chi_{j}^{V_i}$ that are positive definite and radially unbounded w.r.t. the compact set $\mc X^s(\mc C_1)$. Since it holds that $\tilde{\chi}_1 \leq \nu \leq \tilde{\chi}_2$, following the same steps as in ~\cite[p. 98]{H67} there exist $\chi_1 \in \mathscr{K}_{\infty}$ and $\chi_2 \in \mathscr{K}_{\infty}$ such that $\chi_1(\norm{x}_{\mc X^s(\mc C_1)}) \leq \nu(x) \leq \chi_2(\norm{x}_{\mc X^s(\mc C_1)})$. Finally, by the hypothesis of the theorem, the region of attraction $\mathcal{Z}_{\mc X^s(\mc U_1)}$ of $\mc X^s(\mc U_1)$ has measure zero. With $\mc C = \mc X^s(\mc C_1)$, $\mc U = \mc X^s(\mc U_1)$, and $\mc V = \nu$, it follows from Theorem \ref{theorem:AGAS} that \eqref{eq.General.N.dimension} is almost globally asymptotically stable with respect to $\mc X^s(\mc C_1)$. 
\end{IEEEproof}
Theorem \ref{theorem:SPAGAS} requires that the region of attraction $\mc X^s(\mc U_1)$ has measure zero. This can be verified using the following result that relies on the characterization of an unstable hyperbolic fixed point given in Theorem \ref{theorem:ROA} and on Assumption \ref{assumption:derivative}.
\begin{theorem}{\bf(Region of Attraction)} \label{theorem:SPROA}
Suppose that $f_i$ in \eqref{eq.General.N.dimension} is continuously differentiable with linearized dynamics at an equilibrium $x^\star=(x_1,\ldots,x_N)$ given by
\begin{align*}
  \ddt \delta_{x,N} &=  f_{\delta,i}(\delta_{x,1},\ldots,\delta_{x,N}) \coloneqq \sum\nolimits_{j=1}^{N} \tfrac{\partial f_N}{\partial x_j}\big\vert_{x_j=x^\star_j} \delta_{x,j},
\end{align*}
for all $j \in \N_{[1,N]}$, and 
\begin{align*}
 \ddt \delta_{x,i} &=  f_{\delta,i}(\delta_{x,1},\ldots,\delta_{x,i+1}) \coloneqq \sum\nolimits_{j=1}^{i+1} \tfrac{\partial f_i}{\partial x_j}\big\vert_{x_j=x^\star_j} \delta_{x,j}
\end{align*}
for all $i \in \N_{[1,N-1]}$. Suppose that Assumption \ref{assumption:steadystate} and \ref{assumption:derivative} hold for the linearized system (i.e., for $f_{\delta,i}$, $\delta_{x,i}$, $\delta_{y,i}=\delta_{x,i}$ instead of $f_i$, $x_i$, and $y_i$),{ $\mc C_1 = \emptyset$ and $\mc U_1 = \{\mathbbl{0}_{n_1}\}$}, and $V_{\delta,1} = \tfrac{1}{2}\delta^\mathsf{T}_{x,1} P_{\delta,1} \delta_{x,1}$, $V_{\delta,i} =\tfrac{1}{2} \delta^\mathsf{T}_{y,i} P_{\delta,i} \delta_{y,i}$ for all $i \in \N_{[2,N]}$. Moreover, assume that $M$ as in Definition \ref{def.recursive.matrix.M} is positive definite for the linearized dynamics. If there exists $\delta^\prime_{x,1}$ such that $V_{\delta,1}(\delta^\prime_{x,1}) < V_{\delta,1}(\mathbbl{0}_{n_1})$, then the region of attraction $\mc Z_{\{x^\star\},f}$ of $x^\star$ under the nonlinear dynamics \eqref{eq.General.N.dimension} has measure zero.
\end{theorem}
\begin{IEEEproof}
 We define the Lyapunov-like function $\nu_\delta \coloneqq \tfrac{1}{2} \delta^\mathsf{T}_{x,1} P_{\delta,1} \delta_{x,1} + \sum_{i=1}^N  \tfrac{\mu_i}{2}\delta^\mathsf{T}_{y,i} P_{\delta,i} \delta_{y,i}$.
 Following the same steps as in the proof of Theorem \ref{theorem:SPAGAS}, i.e., using Lemma \ref{lemma:decrease.epsilon} (given in Appendix~\ref{app:proofs}) with $\epsilon_i=1$ for all $i \in \N_{[1,N]}$, Lemma \ref{lemma:positive.definiteness}, Assumption \ref{assumption:derivative}, and Proposition \ref{prop:Mdef}, it follows that there exists a a function $\chi_3 \in \mathscr{K}$ such that $\ddt \nu_{\delta} \leq -\chi_3(\norm{\delta_x})$ holds for $\delta_x=(\delta_{x,1},\ldots,\delta_{x,N})$. Moreover, by the hypothesis of the theorem it holds that $V_{\delta,1}(\delta^\prime_{x,1}) < V_{\delta,1}(\mathbbl{0}_{n_1})$ and $V_{\delta,i}(\mathbbl{0}_{n_i})=0$ for all $i \in \N_{[2,N]}$. Letting 
$\delta^\prime_x = (\delta^\prime_{x,1},\mathbbl{0}_{n_2},\ldots,\mathbbl{0}_{n_N})$ it directly follows that $\nu_\delta(\delta^\prime_x) < \nu_\delta(\mathbbl{0}_n)$ and the theorem follows by noting that $\nu_\delta$ satisfies the conditions in Theorem \ref{theorem:ROA}.
\end{IEEEproof}

Note that Theorem \ref{theorem:SPAGAS} and Theorem \ref{theorem:SPROA} provide a way to extend results obtained for the reduced-order system $\ddt x_1 = f^s_1(x_1)$ to the full-order system \eqref{eq.General.N.dimension}. In particular, if there exists a Lyapunov function $V_1$ and Lyapunov-like function $V_{\delta,1}$ for which the conditions of Theorem \ref{theorem:AGAS} and Theorem \ref{theorem:ROA} hold, then Theorem \ref{theorem:SPAGAS} and Theorem \ref{theorem:SPROA} require to find Lyapunov functions $V_i$ and $V_{\delta,i}$ for all $i \in \N_{[2,N]}$ in the error coordinates $y_i$ that satisfy Assumption \ref{assumption:derivative}.

Moreover, Theorem \ref{theorem:SPAGAS} and Theorem \ref{theorem:SPROA} require that $M$ is positive definite. In the next section we exploit the structure of $M$ to provide a recursive sufficient condition for $M$ to be positive definite. This simplified recursive condition will be used to provide analytical stability guarantees for the multi-converter power system considered in Section \ref{sec:model}.

\subsection{Tractable positive definiteness conditions for $M$}\label{sec:poscond}
The following condition is necessary and sufficient for the matrix $M$ in Definition \ref{def.recursive.matrix.M} to be positive definite.
\begin{condition}{\bf(Positivity condition)} \label{condition:stability_conditions}
For all $i\in \mathbb{N}_{[2,N]}$ it holds that $\gamma_i+ \mu_i \boldsymbol{\beta}_i^{\mathsf{T}} M_{i-1}^{-1} \boldsymbol{\beta}^{\phantom{\mathsf{T}}}_i < \alpha_i$, where $M_i$ denotes the $i$-th leading minor of the matrix $M$ as in Definition \ref{def.recursive.matrix.M}.
\end{condition}
Verifying Condition \ref{condition:stability_conditions} requires inverting the matrix $M_i$, i.e., the complexity of the inequalities that need to be verified grows considerably with $N$. In contrast, the following recursive sufficient condition avoids this issue by introducing additional variables $c_i \in \R_{>0}$, that allow us to exploit the structure of the problem. In particular, the variable $c_{i-1}$ lower bounds the smallest eigenvalue of $M_{i-1}$ and can be interpreted as bound on the convergence rate of the reduced-order dynamics \eqref{eq:redordermodel} with $r=i-1$.
\begin{condition}{\bf(Sufficient positivity condition)} \label{condition:sufficient_stability_condition}
For all $i\in \mathbb{N}_{[2,N]}$ it holds that
\begin{align}\label{eq:cipos_cond} 
\alpha_i \beta_{i,i-1} c_{i-1}>\beta_{i-1,i}\boldsymbol{\beta}^\mathsf{T}_i \boldsymbol{\beta}^{\phantom{\mathsf{T}}}_i + \gamma_i \beta_{i,i-1} c_{i-1},
\end{align}
where $c_1 \coloneqq \alpha_1 \in \R_{>0}$, $\alpha_i \in \R_{>0}$, $\gamma_i \in \R$, $\beta_{i,j} \in \R_{>0}$ for all $i\in \mathbb{N}_{[2,N]}$ and $j \in \mathbb{N}_{[1,i-1]}$, and 
\begin{align*} 
c_i \coloneqq \tfrac{1}{2}\left(\alpha_i - \gamma_i + \tfrac{\beta_{i,i-1}}{\beta_{i-1,i}}c_{i-1} - \sqrt{D_i}\right)
\end{align*}
for all $i\in \mathbb{N}_{[2,N]}$ and $D_i \coloneqq (\alpha_i \!-\! \gamma_i \!+\! \tfrac{\beta_{i,i-1}}{\beta_{i-1,i}}c_{i-1})^2 \!+\! 4 (\boldsymbol{\beta}^\mathsf{T}_i \boldsymbol{\beta}_i - (\alpha_i-\gamma_i)\tfrac{\beta_{i,i-1}}{\beta_{i-1,i}}c_{i-1})$.  
\end{condition} 
Condition \ref{condition:sufficient_stability_condition} holds if $y_i$ converges fast enough relative to the reduced-order system \eqref{eq:redordermodel} with $r=i-1$. Next, we show that
\[
 \text{Condition \ref{condition:sufficient_stability_condition}} \Rightarrow \text{Condition \ref{condition:stability_conditions}} \iff M \succ 0.
\]
\begin{proposition}{\bf(Positive definite $\boldsymbol{M}$)} \label{prop:Mdef}
Iff Condition \ref{condition:stability_conditions} holds, then the matrix $M$ is positive definite. Suppose that Condition \ref{condition:sufficient_stability_condition} holds, then {$c_i \in \R_{>0}$ and} $M_i \succeq \mu_i c_i I_i$ holds for all $i\in \mathbb{N}_{[2,N]}$, and the matrix $M$ is positive definite.
\end{proposition}
\begin{IEEEproof}
Using Lemma \ref{lemma:positive.definiteness} (given in Appendix~\ref{app:proofs}) it can be verified that $M$ is positive definite if and only if $M_i\succ0$ (see Definition \ref{def.recursive.matrix.M}) holds for all $i\in\N_{[1,N]}$, i.e., if and only if Condition \ref{condition:stability_conditions} holds. Next, using the partitioning of the matrix $M_i$ from Definition \ref{def.recursive.matrix.M} and applying the Schur complement it can be verified that $M_i \succeq  \mu_i c_i I_i$ holds for all $i \in \N_{[2,N]}$ if and only if $M_{i-1} \succeq \mu_i c_iI_{i-i}$ and
\begin{align} \label{eq:inequality_schur}
  \alpha_i - \gamma_i - c_i - \boldsymbol{\beta}^\mathsf{T}_i (\mu^{-1}_i M_{i-1} - c_i)^{-1} \boldsymbol{\beta}_i \geq 0
\end{align} 
holds for all $i \in \N_{[2,N]}$. The remainder of this proof establishes that \eqref{eq:inequality_schur} and $M_{i-1} \succeq \mu_i c_iI_{i-i}$ hold for all $i \in \N_{[2,N]}$ if Condition \ref{condition:sufficient_stability_condition} holds.
To this end, we first show that if $c_{i-1} \in \R_{>0}$ and Condition \ref{condition:sufficient_stability_condition} hold, then $c_i \in \R_{>0}$ and $\tfrac{\beta_{i,i-1}}{\beta_{i-1,i}}c_{i-1} - c_i \in \R_{>0}$. In particular, $c_i \in \R_{>0}$ holds if and only if  
\begin{align}\label{eq:cipos}
(\alpha_i - \gamma_i + \tfrac{\beta_{i,i-1}}{\beta_{i-1,i}}c_{i-1})^2 > D_i 
\end{align}
holds. Rewriting $D_i>0$ as $D_i = (\alpha_i \!-\! \gamma_i \!-\! \tfrac{\beta_{i,i-1}}{\beta_{i-1,i}}c_{i-1})^2 + 4 \boldsymbol{\beta}^\mathsf{T}_i \boldsymbol{\beta}_i$ it can be verified that \eqref{eq:cipos} is identical to \eqref{eq:cipos_cond}. 
Moreover, by definition of $c_i$ it follows that $\tfrac{\beta_{i,i-1}}{\beta_{i-1,i}}c_{i-1} - 2 c_i = \gamma_i + \sqrt{D_i} - \alpha_i$. Therefore, $\tfrac{\beta_{i,i-1}}{\beta_{i-1,i}}c_{i-1} - c_i \in \R_{>0}$ holds if $c_i+\sqrt{D_i}> \alpha_i - \gamma_i$ holds. Using $c_i \in \R_{>0}$ and $c_{i-1}\in\R_{>0}$ it follows that $\tfrac{\beta_{i,i-1}}{\beta_{i-1,i}}c_{i-1} -c_i \in \R_{>0}$ if \eqref{eq:cipos} holds. This proves the claim. 

Next, we show that $M_{i-1} \succeq \mu_{i-1} c_{i-1}I_{i-1}$, $c_{i-1} \in \R_{>0}$, and Condition \ref{condition:sufficient_stability_condition} guarantee that $M_{i-1} \succeq \mu_i c_iI_{i-1}$ and \eqref{eq:inequality_schur} hold. In particular, $M_{i-1} \succeq \mu_{i-1} c_{i-1}I_{i-1}$ implies that
\begin{align*}
 \mu^{-1}_i M_{i-1} - c_iI_{i-1} \succeq (\frac{\mu_{i-1}}{\mu_i} c_{i-1} - c_i)I_{i-1}. 
\end{align*}

Using $\tfrac{\mu_{i-1}}{\mu_i} c_{i-1} - c_i =\tfrac{\beta_{i,i-1}}{\beta_{i-1,i}} c_{i-1} - c_i > 0,$ it follows that $M_{i-1} - \mu_i c_i I_{i-1} \succeq 0$ and $(\tfrac{\beta_{i,i-1}}{\beta_{i-1,i}}c_{i-1} - c_i)^{-1} I_{i-1} \succeq (\mu^{-1}_i M_{i-1} - c_iI_{i-1})^{-1}$. Hence, \eqref{eq:inequality_schur} holds if 
\begin{align} \label{eq:inequality}
  \alpha_i - \gamma_i - c_i - \boldsymbol{\beta}^\mathsf{T}_i  \big(\tfrac{\beta_{i,i-1}}{\beta_{i-1,i}} c_{i-1} - c_i\big)^{-1} \boldsymbol{\beta}_i \geq 0
\end{align} 
holds. Multiplying \eqref{eq:inequality} by $\tfrac{\beta_{i,i-1}}{\beta_{i-1,i}} c_{i-1} - c_i \in \R_{>0}$ results in the equivalent condition
$\big(\tfrac{\beta_{i,i-1}}{\beta_{i-1,i}} c_{i-1} - c_i\big)\big(\alpha_i - \gamma_i - c_i\big) - \boldsymbol{\beta}^\mathsf{T}_i \boldsymbol{\beta}_i \geq 0$ and it can be verified that $c_i$ defined in Condition \ref{condition:sufficient_stability_condition} is a solution of this inequality. Thus, $M_{i-1} \succeq \mu_{i-1} c_{i-1}I_{i-1}$, $c_{i-1} \in \R_{>0}$, and Condition \ref{condition:sufficient_stability_condition} guarantee that $M_i \succeq \mu_i c_iI_{i}$ holds for all $i \in \N_{[2,N]}$.

Next, noting that $M_1 = \mu_1 c_1 \in \R_{>0}$ ($\mu_1=1$, $c_1=\alpha_1$) it follows by induction that $c_i \in \R_{>0}$ and $M_i \succeq \mu_i c_iI_i$ holds for all $i \in \N_{[2,N]}$. Finally, applying the Sylvester criterion it can be shown that $M_i \succeq \mu_i c_iI_i$ for all $i \in \N_{[1,N]}$ implies that $M \succ 0$.
\end{IEEEproof}

\subsection{Connection to singular perturbation theory}
To illustrate the main differences to results from singular perturbation theory, the following theorem re-states our results in the form typically used in singular perturbation theory.

\begin{theorem}{\bf(Nested Singular Perturbations)}\label{theorem:epsilon}
Consider a compact set $\mc C_1 \subset \R^{n_1}$. Assume that, for all $i\in\N_{[1,N]}$, there exists $\chi^{\scriptscriptstyle{V_i}}_1 \in \mathscr{K}_{\infty}$ and $\chi^{\scriptscriptstyle{V_i}}_2 \in \mathscr{K}_{\infty}$ such that $\chi^{\scriptscriptstyle{V_1}}_1(\norm{x_1}_{\mc C_1}) \leq V_1(x_1) \leq \chi^{\scriptscriptstyle{V_1}}_2(\norm{x_1}_{\mc C_1})$ and 
$\chi^{\scriptscriptstyle{V_i}}_1(\norm{y_i}) \leq V_i(y_i) \leq \chi^{\scriptscriptstyle{V_i}}_2(\norm{y_i})$ for all $i \in\N_{[2,N]}$. {If Assumption \ref{assumption:steadystate} holds and Assumption \ref{assumption:derivative} holds with $\mc U_1 = \emptyset$,} then there exists time constants $0 < \epsilon_N < \epsilon_{N-1} < \ldots < \epsilon_1$ such that the system
	\begin{subequations} \label{eq.epsilon.existence}
	 \begin{align}
		\epsilon_i \ddt x_i &=  f_i (x_1, \hdots, x_{i+1}),  \ \forall i \in \N_{[1,N-1]}, \\
		\epsilon_N \ddt x_N &= f_N (x_1, \hdots, x_{N}), 
	\end{align}
	\end{subequations}	
	is globally asymptotically stable with respect to $\mc X^s(\mc C_1)$.
\end{theorem}
The proof is given in Appendix~\ref{app:proofs}. For $N=2$, $\mc C_1 =\mathbbl{0}_{n_1}$, $\mc U = \emptyset$, and $\psi^\prime_i=0$ for $i \in \N_{[1,2]}$, the conditions of Theorem \ref{theorem:epsilon} and Assumption \ref{assumption:steadystate} and Assumption \ref{assumption:derivative} are identical to the conditions and assumptions in \cite[Th. 11.3]{K02}. Thus, our approach can be interpreted as a generalization of this result to multiple time scales, stability with respect to a compact set, and a less restrictive class of comparison functions $\psi_i$ in Assumption \ref{assumption:derivative}. 
Moreover, we directly bound the convergence rates and interactions between the dynamics on different time scales through Lyapunov functions and do not require both isolating  ``small'' time constants $\epsilon_i$ (see \cite[p. 424]{K02}) and finding suitable Lyapunov functions \cite[Th. 11.3]{K02} to establish stability.

{
\section{Illustrative numerical example}\label{sec:toyexample}
To illustrate the approach presented in the previous section and discuss its implications, we consider the dynamics
\begin{align} \label{eq:toy.example}
\begin{split}
	\!\ddt{} x_1\!&= \!- \tfrac{1}{4} x_1^3-2x_2\eqqcolon f_1(x_1,x_2), \\
	\!\ddt{} x_2 \!&=\! x_1+3x_2 +{\begin{bmatrix}1 & -4\end{bmatrix}} x_3\eqqcolon f_2(x_1,x_2,x_3),\\
	\!\ddt{} x_3 \! &=\!\underbrace{-\! \! \begin{bmatrix} \kappa & 4(1-{\kappa})\\ 0 & 1 \end{bmatrix}\! x_3 \! -\! {\kappa}  \begin{bmatrix} 1+k & 2k+3 \\ 0 & 0 \end{bmatrix}\!\!\begin{bmatrix}                                                                                                            x_1            \\ x_2
\end{bmatrix}}_{\eqqcolon f_3(x_1,x_2,x_3)},
\end{split}
\end{align}
that consist of three nested systems with $n_1=n_2=1$, $n_3=2$, and parameters $k \in \R_{>0}$ and ${\kappa} \in \R_{>0}$.

The steady-state map $\phi_3: \R^{n_1} \times \R^{n_2} \to \R^{n_3}$ of $f_3$ and the vector field $f^s_2 : \R^{n_1} \times \R^{n_2} \to \R^{n_2}$, i.e., $f_2$ with the ``fast'' states $x_3$ restricted to $\phi_3$, are given by
\begin{align*}
\phi_3\! \coloneqq \!-\! \begin{bmatrix} 1+k & 2k+3 \\ 0 & 0 \end{bmatrix}\!\!\begin{bmatrix}                                                                                                            x_1            \\ x_2
\end{bmatrix}
\!, \quad f_2^s\coloneqq-k(x_1 + 2 x_2). 
\end{align*}
Similarly, we obtain the steady-state map $\phi_2: \R^{n_1} \to \R^{n_2}$ of $f^s_2$ and the vector field  $f^s_1 : \R^{n_1} \to \R^{n_1}$, i.e., $f_1$ with the ``faster'' state $x_2$ restricted to $\phi_2$, as $\phi_2 \!\coloneqq \!-\tfrac{1}{2}x_1$ and $f_1^s  \coloneqq -x_1(\tfrac{1}{4} x_1^2-1)$. This results in the ``slow'' reduced-order dynamics of order $n_1$ induced by $f^s_1$, and error coordinates $y_2\coloneqq x_2-\phi_2(x_1)$ and $y_3 \coloneqq x_3-\phi_3(x_1,x_2)$ for the ``fast'' dynamics. Next, we define 
\begin{align*}
 P_3 \coloneqq \begin{bmatrix} \frac{1}{{\kappa}} & -\frac{4}{{\kappa}} \\ -\frac{4}{{\kappa}} & \frac{16}{{\kappa}}+1 \end{bmatrix}
 \end{align*}
and the Lyapunov function candidates $V_1\coloneqq (\frac{1}{4}x^2_1-1)^2$, $V_2\coloneqq \frac{1}{4k}y_2^2$, and $V_3\coloneqq \frac{1}{2} y_3^\mathsf{T}P_3 y_3$, for the reduced-order vector field $f^s_1$, for $x_2$ relative to its steady-state map $\phi_2$, and  $x_3$ relative to its steady-state map $\phi_3$, respectively.

Furthermore, we define the sets $\mc C_1 \coloneqq \{-2,2\}$ and $\mc U_1\coloneqq \{\mathbbl{0} \}$ and the comparison functions $\psi_1 \coloneqq \norm{x_1} \norm{\frac{1}{4} x^2_1-1}$, $\psi_2\coloneqq \norm{y_2}$, $\psi_3\coloneqq\norm{{[1 \; -\!4]} y_3}$, and $\psi_3^\prime\coloneqq \norm{{[0\;\; 1]} y_3}$.%,  with $B^\prime \coloneqq \begin{bmatrix}0 & 1\end{bmatrix}$. 
Next, it can be verified that Assumption 2 holds for $\sigma_1(\norm{x_1}_{\mc C_1 \cup\, \mc U_1}) \coloneqq \tfrac{1}{2}\norm{x_1}_{\mc C_1 \cup\, \mc U_1} \leq \psi_1(x_1)$, and constants $\alpha_1\coloneqq1$, $\alpha_2\coloneqq1$, $\alpha_3\coloneqq\alpha^\prime_3=1$, $\beta_{1,2}\coloneqq 2$, $\beta_{2,3}\coloneqq \frac{1}{2k}$, $\beta_{2,1}=b_{2,1,1}\coloneqq\frac{1}{4k}$, $\gamma_2=b_{2,2,1}\coloneqq -\frac{1}{2k}$, $\beta_{3,1}=b_{3,1,1}\coloneqq \frac{(k+1)}{{\kappa}}$,  $b_{3,2,1}\coloneqq \frac{2(1+k)}{{\kappa}}$, $b_{3,2,2}\coloneqq \frac{2k(2k+3)}{{\kappa}}$, and $\gamma_3=b_{3,3,2}\coloneqq \frac{(2k+3)}{{\kappa}}$. Moreover, note that $\beta_{3,2}=b_{3,2,1}+b_{3,2,2}=\frac{2}{{\kappa}}(2k^2+4k+1)$. 

To apply Theorem \ref{theorem:SPAGAS}, it has to be verified that the matrix $M$,  given in Definition \ref{def.recursive.matrix.M}, is positive definite and the region of attraction of $\mathbbl{0}_{4}$ has measure zero. Condition \ref{condition:stability_conditions} ensures that $M$ is positive definite when 
\begin{align}\label{eq:toy.example.condition1}
k>0, \quad 1>\gamma_3+\frac{(2k+1){\beta^2_{3,1}}+2\beta_{3,1} \beta_{3,2}+{\beta_{3,2}^2}}{2k\beta_{3,2}}.
\end{align}
Condition \ref{condition:sufficient_stability_condition}, is satisfied for $c_2 \coloneqq \frac{\left( 8k+5-\sqrt{64k^2+48k+25} \right)}{16k}$ and 
\begin{align}\label{eq:toy.example.condition2}
k>0, \quad 1>\gamma_3+\frac{{\beta^2_{3,1}}+4{\beta^2_{3,2}}}{8k\beta_{3,2} c_2}.
\end{align}

It remains to show that the region of attraction of $x^\star=\mathbbl{0}_4$ has zero measure. To this end, we define $V_{\delta,1}(\delta_{x,1})=-\frac{1}{2}\delta_{x,1}^2$ and $\psi_1^0(\delta_{x,1})\coloneqq \norm{\delta_{x,1}}$, and note that $f_2^s$, $f_3^s$ are linear and corresponding Lyapunov function candidates are quadratic. Using $\psi_2$, $\psi_3$, $\psi_3^\prime$ and the same constants as above, it can be verified that Assumption 2 holds for the linearization of \eqref{eq:toy.example}. Thus, when either \eqref{eq:toy.example.condition1} or \eqref{eq:toy.example.condition2} holds, it follows that $\ddt{}\nu_{\delta_x}<0$ for all $\delta_x\neq \mathbbl{0}_{4}$. Lastly, for $\delta_{x}^\prime \coloneqq (1,\mathbbl{0}_3)$ it holds that $\nu_{\delta_x}(\delta_{x}^\prime) < \nu_{\delta_x}(\mathbbl{0}_4)$. Hence, it follows from Theorem \ref{theorem:SPROA} that the region of attraction of $\mathbbl{0}_4$ has measure zero and Theorem 3 guarantees that the system \eqref{eq:toy.example} is almost globally asymptotically stable with respect to $\mc X^s(\mc C_1)$.

To illustrate the stability conditions we define the set of parameters 
\begin{align*}
 \mc P_4 \! = \! \left\{ \! \begin{bmatrix} {\kappa} \\  k\end{bmatrix} \!\in\! \R^2_{>0} \left\vert {\kappa} \!> \!8 +\! 4k\!+ \!\frac{2k^3+21k+36k+9}{4k(2k^2+4k+1)}\! \right.\right\}
\end{align*}
that satisfy Condition \ref{condition:stability_conditions} (i.e., \eqref{eq:toy.example.condition1}) and the set 
\begin{align*}
 \mc P_3 \!= \! \left\{ \! \begin{bmatrix} {\kappa} \\  k\end{bmatrix} \!\in\! \R^2_{>0} \left\vert {\kappa} \!> \!3\! + \! 2k\!+\! \frac{(1\!+\!k)^2+16(2k^2\!+\!4k\!+\!1)^2}{16k(2k^2\!+\!4k\!+\!1)c_2} \right.\right\}
\end{align*}
of parameters that satisfy Condition \ref{condition:sufficient_stability_condition} (i.e., \eqref{eq:toy.example.condition2}).
Figure \ref{fig:stabilityregion} shows different regions of the parameter space. For $({\kappa},k) \in \mc P_4$ Condition \ref{condition:stability_conditions} holds and for $({\kappa},k) \in \mc P_3$ Condition \ref{condition:sufficient_stability_condition} holds and \eqref{eq:toy.example} is almost globally asymptotically stable with respect to $\mc C_1$. Moreover, for $({\kappa},k) \in \mc P_6$ the linearization at $x_1 = \pm 2$ is unstable, i.e., the system cannot be almost globally asymptotically stable with respect to $\mc C_1$. For $({\kappa},k) \in \mc P_5$ we can neither certify stability nor instability. 

The analysis above can be repeated for $\psi_3=\norm{y_3}$ and $\psi^\prime_3=0$, i.e., when not separating the subspace of $x_3$ that directly acts on $x_2$. In this case, Condition \ref{condition:stability_conditions} corresponds to $({\kappa},k) \in \mc P_2$ and Condition \ref{condition:sufficient_stability_condition} corresponds to $({\kappa},k) \in \mc P_1$. We note that $\mc P_1 \subseteq \mc P_2$, $\mc P_2 \subseteq \mc P_3$, and $\mc P_3 \subseteq \mc P_4$. In other words, comparing either $\mc P_1$ and $\mc P_2$ or $\mc P_3$ and $\mc P_4$ we observe that Condition \ref{condition:sufficient_stability_condition} is more conservative than Condition \ref{condition:stability_conditions}. Moreover, comparing either $\mc P_1$ and $\mc P_3$ or $\mc P_2$ and $\mc P_4$, it can be seen that exploiting the degree of freedom introduced by $\psi^\prime_i$ in Assumption \ref{assumption:derivative} results in less conservative bounds. {While the precise stability boundary in the parameter space is not known, we observe that the condition $({\kappa},k) \in \mc P_4$ is only mildly conservative (i.e., at most a factor of 10).}

\begin{figure}[t!!!]
\centering
\vspace{0em}
\includegraphics[width=1\columnwidth]{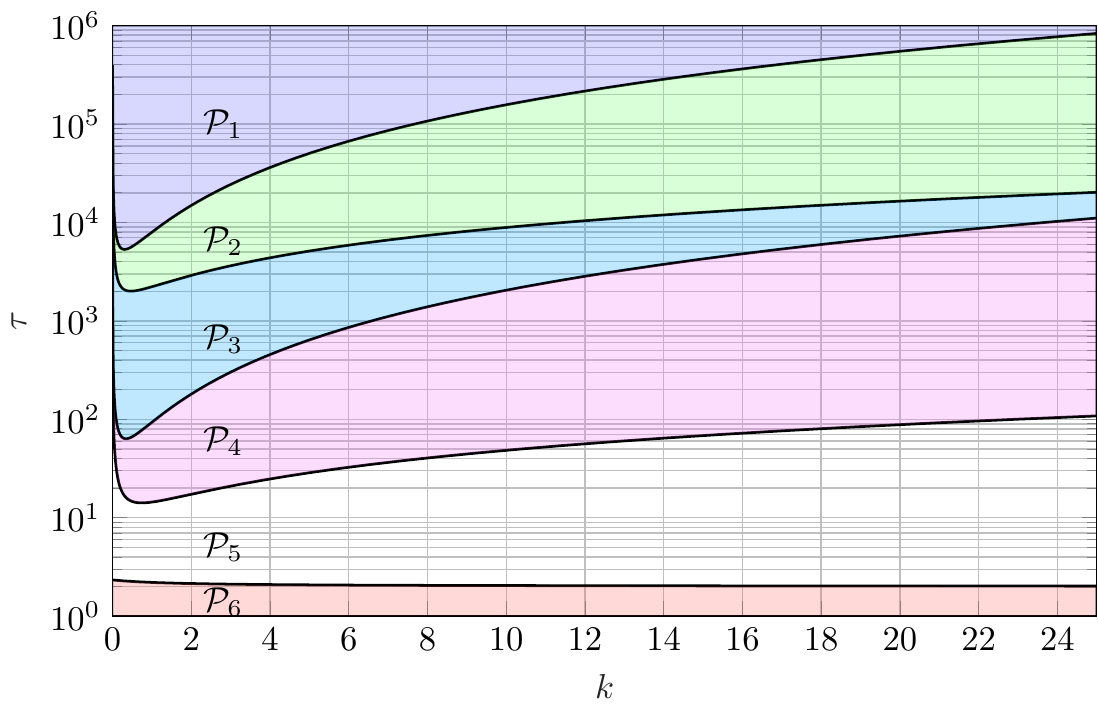}
\vspace{-2em}
\caption{{Stability regions in the parameter space (note that $\mc P_1 \subseteq \mc P_2 \subseteq \mc P_3 \subseteq \mc P_4$). For   $\psi_3=\norm{{[1\; -\!4]} y_3}$, and $\psi_3^\prime=\norm{{[0\;\; 1]} y_3}$ we obtain the stability conditions $({\kappa},k) \in \mc P_4$ (Condition \ref{condition:stability_conditions}) and $({\kappa},k) \in \mc P_3$ (Condition \ref{condition:sufficient_stability_condition}). For $\psi^\prime_3=0$ and $\psi_3=\norm{y_{{3}}}$ we obtain the stability conditions $({\kappa},k) \in \mc P_2$ (Condition \ref{condition:stability_conditions}) and $({\kappa},k) \in \mc P_1$ (Condition \ref{condition:sufficient_stability_condition}). For $({\kappa},k) \in \mc P_5$ we can neither certify stability nor instability. \label{fig:stabilityregion}}}
\end{figure}
}

\section{Multi-converter power system and control objectives}\label{sec:model}

\subsection{Model of a converter based power system}
\begin{figure}[b!!!]
	\begin{center}
		\vspace*{-0.25cm}
		\includegraphics[width=0.95\linewidth]{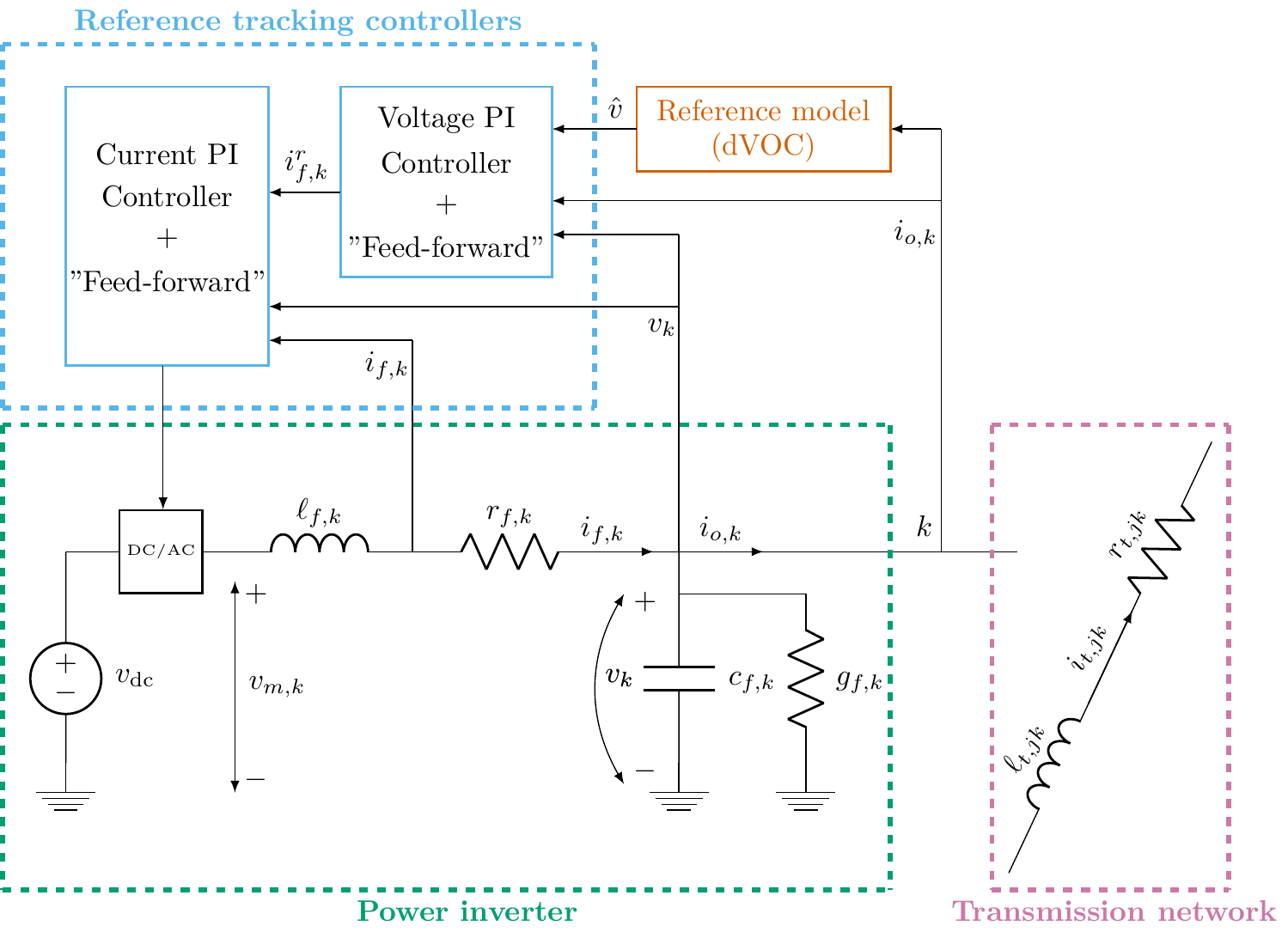}
	\end{center}
\vspace*{-0.25cm}
\caption{Model of a DC/AC converter with reference model, reference tracking controller, and transmission network.}
\label{fig:Model}
\end{figure}

The three-phase power system considered in this article consists of ${N_c}$ DC/AC converters, interconnected through ${N_t}$ resistive-inductive lines. {All electrical signals are assumed to be balanced and, for notational convinience, we work in a global $dq$-coordinate frame with reference angle $\theta_r = \omega_0 t$ that rotates with the nominal frequency $\omega_0$ \cite[Ch. 6]{PWS-MAP:98}. We emphasize that \textit{(i)} we model electrical signals (i.e., voltages in currents) in two-dimensional rectangular coordinates and not in polar coordinates, and \textit{(ii)} the implementation of the controls in stationary $abc$ or $\alpha\beta$-coordinates \cite{clarke1943circuit} are independent of the reference angle $\theta_r$ allowing for a decentralized implementation (please see \cite{colombino2017global2, line_dynamics} for further details).

We consider an averaged model of a two-level voltage source converter that is interfaced to the transmission network via a $RLC$ filter. The overall setup is depicted in Figure \ref{fig:Model}. The network topology is modeled by the oriented incidence matrix $B \in \{-1, 0, 1\}^{N_c \times N_t}$ of the connected, undirected, weighted graph $\mc G$ with $N_c$ nodes corresponding to the converters and the set of edges $\mc E \subseteq \N_{[1,N_c]}\times \N_{[1,N_c]}$ corresponding to the transmission lines \cite[Section II.A]{line_dynamics} with inductance $\ell_{t,jk} \in \R_{> 0}$ and resistance $r_{t,jk} \in \R_{> 0}$ for all $(j,k) \in \mc E$.  For brevity of the notation, we associate an index $l \in \{1,\hdots N_t \}$ to each edge. Using the extended incidence matrix  $\mc B \coloneqq B \otimes I_2$, the dynamics of the $N_t$ transmission line currents ${i}_t \coloneqq (i_{t,1}, \hdots i_{t,{N_t}}) \in \R^{2N_t}$ in the rotating frame are given by (see \cite[Sec. II]{line_dynamics})
\begin{align} \label{eq:transmission_lines_dynamics_rotating}
L_T \ddt{} i_t = -Z_T i_t + \mc B^\mathsf{T} v,
\end{align} 
with $J_k \coloneqq I_{k} \otimes J$, and impedance, resistance, and inductance matrices $Z_T{\coloneqq} R_T+\omega_0 J_{N_t} L_T$, $R_T \coloneqq \diag (\{r_{t,l} {I_2}\}_{l=1}^{N_t})$, and $L_T \coloneqq \diag (\{\ell_{t,l} {I_2}\}_{l=1}^{N_t})$. {Moreover, we define the network admittance matrix $\mc \Yn \coloneqq \mc B Z^{-1}_T \mc B^\mathsf{T}$ and the symmetric Laplacian and extended Laplacian matrices are defined by $\mc L \coloneqq L \otimes I_2$ and $\mc L \coloneqq \mc R(\kappa) \mc \Yn$ (see \cite[Sec. II]{line_dynamics} for further details).}

Next, we consider an averaged model of a two-level voltage source converter that modulates a DC voltage $v_{\text{DC},k}$ into an AC voltage $v_{m,k}$ (see \cite{CGD17, JAD18}). Additionally, we assume that the DC voltage $v_{\text{DC},k}$ is regulated to be constant by the source feeding the converter. Using the vector of converter terminal voltages $v \coloneqq (v_1, \hdots v_{N_c}) \in \R^{2N_c}$, filter currents $i_f\coloneqq (i_{f,1}, \hdots i_{f,{N_c}}) \in \R^{2N_c}$, and modulated AC voltages $v_m\coloneqq (v_{m,1}, \hdots v_{m,{N_c}}) \in \R^{2N_c}$, and the vector $i_o = (i_{o,1},\ldots,i_{o,N_c}) = \mc B i_t$ of converter current injections, the converters' open loop dynamics are given by 
\begin{subequations} \label{eq:filter_dynamics_rotating_matrix}
	\begin{align} 
	L_f \ddt{} i_f &= -Z_f i_f - v + v_m, \label{eq:if_rf.matrix}\\
	C_f \ddt{} v &= - Y_f v - i_o + i_f \label{eq:v_rf.matrix},
	\end{align}
\end{subequations}
where $Z_f = R_f + J_{N_c} \omega_0 L_f$ and $Y_f = G_f + J_{N_c} \omega_0 C_f$ model the filter impedance and admittance, and the matrices $R_f \coloneqq \diag (\{r_{f,k} {I_2}\}_{k=1}^{N_c})$, $L_f \coloneqq \diag (\{\ell_{f,k}{I_2}\}_{k=1}^{N_c})$, $C_f \coloneqq \diag (\{c_{f,k}{I_2}\}_{k=1}^{N_c})$, and $G_f \coloneqq \diag (\{g_{f,k}{I_2}\}_{k=1}^{N_c})$ model the filter resistance, inductance, capacitance, and conductance.
}
\subsection{Operating point and control objectives} \label{sec:control_objective}%
The main control objective is to steer the system to an operating point that is commonly defined by the nominal frequency $\omega_0$, set-points $v^\star_k \in \R_{\geq 0}$ for the terminal voltage magnitude $\norm{v_k}$, and set-points $p^\star_k \in \R$ and $q^\star_k \in \R$ for the active and reactive power injections $p_k \coloneqq v_{k}^\mathsf{T}  i_{o,k} \in \mathbb{R}$ and $q_k \coloneqq v_{k}^\mathsf{T}  J i_{o,k} \in \mathbb{R}$ (see \cite[Def. 2]{line_dynamics}). Note that the terminal voltage in stationary coordinates (i.e., $abc$ or $\alpha\beta$) oscillates with the nominal frequency and constant amplitude if $\ddt v_k = 0$. We emphasize that stability of the nominal operating point does not imply that the converters operate with constant power control. Instead, the terminal voltages are adjusted to stabilize the power system at the given operating point. This operating point cannot be prescribed arbitrarily, but needs to satisfy the power-flow equations \cite{PWS-MAP:98}.
\begin{condition}{\bf(Consistent operating point)} \label{cond.consistent}
	The set-points $p_k^\star \in \mathbb{R}$, $q_k^\star \in \mathbb{R}$, $v_k^\star \in \mathbb{R}{>0}$ for active power, reactive power, and voltage magnitude respectively, are consistent with the power flow equations with line conductance and susceptance $g_{t,jk}=\tfrac{r_{t,jk}}{r^2_{t,jk}+\omega^2_0\ell^2_{t,jk}}$ and $b_{t,jk}=\tfrac{\omega_0\ell_{t,jk}}{r^2_{t,jk}+\omega^2_0\ell^2_{t,jk}}$, i.e., for all $(j,k) \in \mc E$ there exist relative voltage phase angles ${\theta^\star_{jk}} \in [-\pi,\pi]$ such that branch powers
	\begin{align*}
	p^\star_{jk} &\coloneqq v^{\star2}_k g_{t,jk} - v^{\star}_j v^{\star}_k
		( g_{t,jk}\cos(\theta_{jk}^\star )+b_{t,jk}\sin(\theta_{jk}^\star ),\\
	q^\star_{jk} &\coloneqq v^{\star2}_k b_{t,jk} - v_k^{\star}v_j^{\star}( b_{t,jk} \cos(\theta_{jk}^\star )-g_{t,jk} \sin(\theta_{jk}^\star ),
	\end{align*}
	at the nominal operating point satisfy $p^\star_k = \sum_{(j,k)\in\mc E} p^\star_{jk}$ and $q^\star_k = \sum_{(j,k)\in\mc E} q^\star_{jk}$ for all buses $k \in \mathbb{N}_{[1,N]}$.
\end{condition}
Note that Condition \ref{cond.consistent} uses the standard power flow equations in polar coordinates and $\theta_{jk}^{\star}$ denotes the angle difference between the terminal voltages of the $j^\text{th}$ and $k^\text{th}$ converter at the nominal operating point. In contrast, the state $v_k$ corresponding to the converter terminal voltage is defined in rectangular $dq$-coordinates. Thus, at the nominal operating point, the converter terminal voltage needs to satisfy $v_j/v^\star_j = \mc R(\theta^\star_{jk}) v_k/v^\star_k$ for all $(j,k)\in\mc E$. We emphasize that the controllers developed in the next Section will not use the relative angles $\theta^\star_{jk}$, but only require $p^\star_k$, $q^\star_k$, and $v^\star_k$ to define the nominal operating point.

In contrast, for the purpose of stability analysis, it is more convinient to convert the local operating point specification $p^\star_k = v^\mathsf{T}_k i_{o,k}$, $q^\star_k = v^\mathsf{T}_k J i_{o,k}$, $\norm{v_k}=v^\star_k$ into the equivalent non-local set-based specification.
\begin{subequations}\label{eq:formspecs}
\begin{align}
\mc S &\coloneqq \left\{v \in \R^{2N_c} \left\vert\; \frac{v_k}{v^\star_k}=  R(\theta_{k1}^\star)\,\frac{v_1}{v^\star_1},\; \forall k \in \N_{[2,N_c]} \right.\right\},\\
\mc A &\coloneqq \left\{v \in \R^{2N_c} \left\vert\; \norm{v_k} =  v_k^\star,\; \forall k \in \N_{[1,N_c]} \right.\right\}.
\end{align} 
\end{subequations}
The intersection $\mc S \cap \mc A$ of the sets $\mc S$ (voltage phase difference) and $\mc A$ (voltage magnitude) encodes the nominal operating point for the reference voltages \cite[Sec. II-C]{line_dynamics}. Moreover, if the voltage dynamics converge to their steady-state map, i.e., $x_v = \phi_v(\hat{v}) = (\hat{v},\mathbbl{0}_{2 N_c})$, it follows that $v \in \mc S \cap \mc A$ and the converter terminal voltages $v$ are at the prescribed operating point. Finally, the full steady-state specifications in the rotating frame can be expressed as $x \in \mc X^s(\mc S \cap \mc A)$.

\section{Dispatchable virtual oscillator control with cascaded inner control loops}\label{sec:control_design}
In this section, we propose a control law that admits a fully decentralized implementation and consists of two cascaded inner loops that track reference voltage. Typically nonlinear stability analysis of multi-converter systems is severely complicated by standard cascaded inner loops that are implemented in local rotating reference frames for every converter \cite{PPG07}. In contrast, we propose a cascaded two-degree of freedom inner control structure that can be directly implemented in stationary (i.e., non-rotating) coordinates and, under suitable assumptions, preserves the stability guarantees obtained for dispatchable virtual oscillator control obtained under the assumption that the terminal voltage $v_k$ can be controlled directly \cite{colombino2017global2,line_dynamics}. After presenting the control law we provide the steady-state maps and reduced-order models required to apply the results developed in Section \ref{section:SPT.extension}.

\subsection{Dispatchable virtual oscillator control as reference model}
We require the following assumption that is commonly made in the stability analysis of AC power systems.
\begin{assumption}{\bf(Uniform inductance-resistance ratio)}\label{ass.constant.ratio}
	The inductance to resistance ratio of every transmission line is constant, i.e., for all $(j,k) \in\mc E$ it holds that $\frac{\ell_{jk}}{r_{jk}} = \rho \in \R_{>0}$.
\end{assumption}
This assumption typically holds for transmission lines on the same voltage level and simplifies the analysis while preserving the main salient features of the system. Additionally, we define the angle $\kappa \coloneqq \tan^{-1}(\omega_0 \rho)$. Next, consider the vector $\hat{v}=(\hat{v}_1,\ldots,\hat{v}_{N_c}) \in \R^{2N_c}$ of voltage references and the dynamics
\begin{align} \label{eq.closed.loop.r.f.hatv}
\ddt \hat{v} = \eta\big(\mc K \hat{v} -\mc R(\kappa) \mc B i_t +\eta_a\Phi(\hat{v}) \hat{v}\big), 
\end{align}
where $\eta\in \R_{>0}$ and $\eta_a\in R_{>0}$ are control gains, and \begin{align*}
\mc K &\coloneqq \diag(\{K_k\}_{k=1}^{N_c}),
&\Phi(\hat{v})&\coloneqq \diag(\{\Phi(\hat{v}_k)I_2\}_{k=1}^{N_c}),
\end{align*} with $K_k$ and $\Phi_k$ defined as
\begin{align}
K_k \!\coloneqq\! \frac{1}{{v_k^\star}^2} \mc R(\kappa)\! \begin{bmatrix}
p_k^{\star} & q_k^{\star}\\
-q_k^{\star} & p_k^{\star}
\end{bmatrix}\!, \; \Phi_k(\hat{v}_k)\coloneqq 1 - \frac{\norm{\hat{v}_k}^2}{v_k^{\star 2}} .\label{eq:Kphi}
\end{align}
The multi-converter system (without internal converter dynamics) is almost globally asymptotically stable with respect to the specifications given in \ref{sec:control_objective}, if $v=\hat{v}$, i.e., if the converter terminal voltages evolve according to the reference dynamics \eqref{eq.closed.loop.r.f.hatv} (see \cite{line_dynamics}). Therefore, we will use $\hat{v}$ as a reference for an underlying voltage controller. Note that $\ddt v =0$ implies that $ v$ rotates with the nominal frequency $\omega_0$. Moreover, $\Phi_k( v_k) v_k$ can be interpreted as a voltage regulator, i.e.,depending on the sign of the normalized quadratic voltage error $\Phi_k( v_k)$ the voltage vector $ v_k$ is scaled up or down. Finally, the term $K_k v_k - R(\kappa)i_{o,k}$ can be interpreted in terms of  power set-point tracking or phase synchronization (see \cite[Sec. II-D]{line_dynamics}).

\subsection{Two-degree-of-freedom reference tracking controllers} 
We propose two local (i.e., decentralized) cascaded two-degree-of-freedom PI controllers for the filter currents and terminal voltage. An outer loop provides a reference $i_{f}^r$ for the filter currents $i_{f}$ and ensures that the terminal voltages $v$ are tracking the reference $\hat v$. Furthermore, an inner loop computes the control signal ${v}_{m}$ so that the filter current $i_{f}$ tracks the reference $i_{f}^r$, provided by the outer loop.

We first assume that the filter currents $i_{f}$ can be used to control the voltage dynamics \eqref{eq:v_rf.matrix} and propose the two-degree of freedom voltage PI controller
\begin{subequations}\label{eq:voltage.controller.single}
	\begin{align}
	 \ddt{}\zeta_{v}  \!\! &\coloneqq v- \hat{v}, \\
	{i}_{f}^r &\coloneqq Y_{f} v + i_{o} - K_{p,v}(v - \hat{v}) - K_{i,v} \zeta_{v},
	\end{align}
\end{subequations}
where the term $Y_{f} v$ compensates the filter admittance losses at the nominal operating point, $\zeta_{v} \in \R^{2N_c}$ denotes the voltage controller integrator states, and $K_{p,v} \coloneqq  \diag(\{K_{p,v_{k}} I_2\}_{k=1}^{N_c})$ and $K_{i,v} \!  \coloneqq \!  \diag(\{K_{i,v_{k}} I_2\}_{k=1}^{N_c})$ collect the  proportional and integral gains $K_{p,v_k} \in \R_{>0}$ and $K_{i,v_k} \in \R_{>0}$ of the individual converters. We stress that each voltage controller obtains $v$, $i_o$ and $i_{f}$ from the local measurements. However, because $i_{f}$ is not the system's input, the controller \eqref{eq:voltage.controller.single} cannot be applied directly. Therefore, we use another two-degree of freedom current PI controller 
\begin{subequations} \label{eq:filter.current.controller.single}
	\begin{align}
	\ddt{}\zeta_{f} &\coloneqq i_{f}-i^r_{f} \\
	v_{m} & \coloneqq Z_{f}  i_{f} + v - K_{p,f}(i_{f} -i_{f}^r)- K_{i,f} \zeta_{f},
	\end{align}
\end{subequations}
that tracks the reference signals $i_{f}^r$ by acting on the control inputs $v_{m}$. The term $Z_{f}  i_{f}$ compensates the filter impedance losses at the nominal operating point, $\zeta_{f}\in \R^{2N_c}$ denotes the integrator states, and $K_{p,f} \coloneqq \diag(\{K_{p,f_k} I_2 \}_{k=1}^{N_c})$ and $K_{i,f} \coloneqq \diag(\{K_{i,f_k} I_2 \}_{k=1}^{N_c})$ collect the proportional and integral gains $K_{p,f_k} \in \R_{>0}$ and $K_{i,f_k} \in \R_{>0}$ of the converters. 

\begin{remark}{\bf(Decentralized control implementation)}
The controllers \eqref{eq.closed.loop.r.f.hatv}, \eqref{eq:voltage.controller.single}, and  \eqref{eq:filter.current.controller.single} are independent of the $dq$-frame reference angle $\theta_r$ when implemented in stationary $\alpha\beta$-coordinates. Moreover, we do not require clock synchronization because the implementation in $\alpha\beta$-coordinates does not require reference angles (i.e., $\theta^\star_{jk}$). Instead, we rely on droop like properties for angle synchronization \cite[Sec. III-D]{colombino2017global2}.      
     \end{remark}

\subsection{Closed-loop dynamics and steady-state maps}
Note that the closed-loop multi-converter system 
\begin{subequations}\label{eq:closed.loop.r.f.dynamcics}
	\begin{align}
	\ddt \hat{v} &= \underbrace{\eta\big(\mc K \hat{v} -\mc R(\kappa) \mc B i_t +\eta_a\Phi(\hat{v}) \hat{v}\big)}_{\eqqcolon f_{\hat{v}}(\hat{v},i_t)}, \label{eq.closed.loop.hatv}  \\
	\ddt{i_t} &= \underbrace{ L^{-1}_T  (-Z_T i_t + \mc B^\mathsf{T} v)}_{\eqqcolon f_t(i_t,x_v)}, \label{eq.closed.loop.it} \\
	  \ddt{}  \underbrace{\begin{bmatrix}
		v \\ \zeta_v
		\end{bmatrix}}_{\eqqcolon x_v}& = \underbrace{\begin{bmatrix} C_f^{-1} (- Y_f v - \mc B i_t + i_f) \\ v-\hat{v}\end{bmatrix}}_{\eqqcolon f_v(\hat{v},i_t,x_v,x_f)}, \label{eq.closed.loop.xv} \\
	 \ddt{}  \underbrace{\begin{bmatrix}
		i_f \\ \zeta_{f}
		\end{bmatrix}}_{\eqqcolon x_f} &=\underbrace{\begin{bmatrix} -L_f^{-1} ( -K_{p,f} (i_f-i_f^r) -K_{i,f} \zeta_{f}) \\ i_f-i_f^r\end{bmatrix}}_{\eqqcolon f_f(\hat{v},i_t,x_v,x_f)}, \label{eq.closed.loop.r.f.xif}
	\end{align}
\end{subequations}
is in the form of the nested dynamical system \eqref{eq.General.N.dimension} with $N\!=\!4$, $(x_1,x_2,x_3,x_4) \!=\! (\hat{v},i_t,x_v,x_f)$, $(f_1,f_2,f_3,f_4) \!=\! (f_{\hat{v}},f_t,f_v,f_f)$, and the steady-state maps $(\phi_1,\phi_2,\phi_3) \coloneqq (\phi_f,\phi_v,\phi_t)$ defined in Assumption \ref{assumption:steadystate} are given by
\begin{subequations} \label{eq:ss.maps}
	\begin{align} 
	\phi_{f} (\hat{v}, i_t,x_v)  &\coloneqq (i_f^r(\hat{v},i_t,x_v), \mathbbl{0}_{2N_c}), \label{eq:ss.map.xif}\\
	\phi_{v}(\hat{v}) &\coloneqq (\hat{v}, \mathbbl{0}_{2N_c}), \label{eq:ss.map.xv}\\
	\phi_t(\hat{v}) &\coloneqq Z_t^{-1}\mc B^{\mathsf{T} } \hat{v} \label{eq:ss.map.it}. 	
	\end{align} 
\end{subequations}
Moreover, the functions $f_v^s$, $f_t^s$, and $f_{\hat{v}}^s$ are given by
\begin{subequations} \label{eq:f.maps}
	\begin{align} 
	f_v^s(\hat{v},x_v) &\coloneqq \begin{bmatrix}
	- C_f^{-1} (K_{p,v}(v-\hat{v})+K_{i,v} \zeta_{v}) \\ v-\hat{v}
	\end{bmatrix}, \label{eq:ss.f.xv}\\
	f_t^s(\hat{v},i_t) &\coloneqq L_T^{-1}(-Z_t^{-1} i_t +\mc B^{\mathsf{T} } \hat{v} )\label{eq:ss.f.it}, \\
	f_{\hat{v}}^s(\hat{v}) &\coloneqq \eta\big(\mc K \hat{v} -\mc R(\kappa)  \mc \Yn \hat{v} +\eta_a\Phi(\hat{v}) \hat{v}\big). \label{eq:ss.f.vhat}
	\end{align} 
\end{subequations}
Note that \eqref{eq:ss.f.xv} is obtained by assuming that the current control \eqref{eq:filter.current.controller.single} perfectly tracks its reference. Moreover, \eqref{eq:ss.f.it} is obtained by additionally assuming that the voltage control \eqref{eq:voltage.controller.single} perfectly tracks its reference. Finally, \eqref{eq:ss.f.vhat} is the reduced-order model of the closed-loop system obtained by assuming that the inner controls perfectly track their reference and that the transmission line dynamics are negligible. This is the setup considered in \cite{colombino2017global2}. Note that $i_t = \phi_t(v)$ is the quasi-steady-state network model that is commonly used for stability analysis of conventional power systems. This approximation is typically justified due to the pronounced time-scale separation between the dynamics of the transmission lines and the dynamics of synchronous machines. However, for converter-based power systems the electromagnetic transients of the lines have a significant influence on the stability boundaries, and the approximation is no longer valid \cite{vorobev2017high,line_dynamics,Guerrero2015,MSV+19}.

\subsection{Main result}\label{section:Stability.Analysis}
{In the following, we provide conditions under which the multi-converter system \eqref{eq:closed.loop.r.f.dynamcics} is almost globally asymptotically stable with respect to $\mc X^s(\mc S \cap \mc A)$ (i.e., $\hat{v} \in \mc S \cap \mc A$, $v=\hat{v}$, $i_t=\phi_t(v)$, $i_f=i^r_f$). This implies almost global asymptotic stability of the operating point specified by $v \in \mc S \cap \mc A$ in \eqref{eq:formspecs}.
\begin{condition}{\bf(Stability Condition)}\cite[Prop. 2]{line_dynamics}\label{cond.stab}
 The set-points $p^\star_k$, $q^\star_k$, $v_k^\star$, and the branch powers $p^\star_{jk}$, $q^\star_{jk}$ at the nominal operating point satisfy Condition~\ref{cond.consistent}, $|\theta^\star_{jk}| \leq \tfrac{\pi}{2}$ holds for all $(j,k) \in \mc N \times \mc N$, and for all $k \in \mc N$, the network load margin $c_L \in \mathbb{R}_{>0}$, the network parameters $\norm{Y_{jk}}$, and the gains $\eta \in \mathbb{R}_{>0}$, $\alpha \in \mathbb{R}_{>0}$ satisfy
 \begin{align*}
  \sum_{j : (j,k) \in \mc E} \frac{\cos(\kappa)}{v^{\star2}_k} \left|p^\star_{jk}\right|\!+\!\frac{\sin(\kappa)}{v^{\star2}_k} \left| q^\star_{jk} \right|\!+\! \eta_\alpha \leq 
\frac{v^{\star2}_{\min}}{2 v^{\star2}_{\max}}\lambda_2( L)-c_L,\\
\eta < \frac{c_L}{2 \rho d_{\max} (c_L  + 5 \max_{k \in \mc N} \sqrt{p^{\star2}_k+q^{\star2}_k}  {v^\star_k}^{-2} + 10 d_{\max})},
 \end{align*}
where $d_{\max} \!\coloneqq\! \max_{k \in \mc N} \sum\nolimits_{j:(j,k) \in \mc E} \norm{Y_{jk}}$ is the maximum weighted node degree of the transmission network graph, and $\lambda_2(L)$ is the second smallest eigenvalue of the graph Laplacian $L$.
\end{condition}
Broadly speaking, the first inequality in Condition \ref{cond.stab} requires that the network is not too heavily loaded. The margin $c_L$ links the network loading to the second condition that enforces a sufficient time-scale separation between the dVOC reference model and the network dynamics. Note that the bound on $\eta$ decreases with higher network loading (i.e., larger $c_L$), increased power injection (i.e., larger apparent power $\sqrt{p^{\star2}_k+q^{\star2}_k}$) and connectivity of individual nodes (i.e., larger $d_{\max}$).}
\begin{condition} \label{cond:voltage.loop.gains}
Consider the control gain $\eta$ and constant $c_L$ that satisfy Condition \ref{cond.stab} and the stability margin $c_2$ defined in Condition \ref{condition:sufficient_stability_condition} (with $c_1=\alpha_1$, $\alpha_2=1$, $\gamma_2=\eta  \rho  \norm{\mc \Yn}$, $\beta_{1,2} = \norm{\mc K - \mc L}^{-1}$, and $\beta_{2,1} = \rho \norm{\mc \Yn}$). The control gains $K_{p,v_k} \in \R_{>0}$, $K_{i,v_k} > c_{f,k}$ of the voltage PI controller satisfy 
\begin{align*}%\label{eq:voltagegainbound}
	\frac{1+ \max\limits_{k \in \N_{[1,N_c]}}  \tfrac{K_{i,v_k}}{K_{p,v_k}}}{\min\limits_{k \in \N_{[1,N_c]}} \tfrac{K_{i,v_k}}{c_{f,k}} -1 } < \frac{4 \eta c_2}{\norm{\mc B R_T^{-1} \mc B^{\mathsf{T}}} (1+4 \eta^2)}.
\end{align*}
\end{condition}
\begin{condition} \label{cond:current.loop.gains}
Consider the control gains $\eta$, $K_{p,v_k}$, and $K_{i,v_k}$, and stability margins $c_1$, $c_2$ such that Condition \ref{cond:voltage.loop.gains} is satisfied. Given the stability margin $c_3$ defined in Condition \ref{condition:stability_conditions} (with $\alpha_3 = 1-\max_{k \in \N_{[1,N_c]}} \frac{c_{f,k}}{K_{i,v_k}}$, $\gamma_3=0$, $\beta_{3,2}=\eta \beta_{3,1}$, $\beta_{3,1}=\max_{k \in \N_{[1,N_c]}} \tfrac{c_{f,k}}{K_{i,v_k}}+\tfrac{c_{f,k}}{K_{p,v_k}}$, and $\beta_{2,3}=\norm{\mc B R^{-1}_T \mc B^\mathsf{T}}$), the control gains $K_{p,f_k} \in \R_{>0}$, $K_{i,f_k} > \ell_{f,k}$ of the current PI controller satisfy
\begin{align*}%\label{eq:currentgainbound}
\frac{1+ \max\limits_{k \in \N_{[1,N_c]}}  \tfrac{K_{i,f_k}}{K_{p,f_k}}}{\min\limits_{k \in \N_{[1,N_c]}} {\tfrac{K_{i,f_k}}{\ell_{f,k}}} -1} < \frac{4 c_3}{\tfrac{\beta_{34}}{\tilde{\beta}_{43}}(\tilde{\beta}^2_{41}+\tilde{\beta}^2_{42}+4 \tilde{\beta}_{43}^2)+c_3 \tilde{\gamma}_4},
\end{align*}
where $\beta_{34}=\max_{k \in \N_{[1,N_c]}} \tfrac{1}{K_{i,v_k}}+\tfrac{1}{K_{p,v_k}}$, $\tilde{\beta}_{41}=\max_{k \in \N_{[1,N_c]}} K_{p,v_k}$, $\tilde{\beta}_{42}=\tfrac{\omega_0}{\sin(\kappa)}+\eta \max_{k \in \N_{[1,N_c]}} K_{p,v_k}$, $\tilde{\beta}_{43}=\norm{Y_f-K_{p,v}}\max_{k \in \N_{[1,N_c]}}(\tfrac{K_{p,v_k}}{c_{f,k}}+\tfrac{K_{i,v_k}}{c_{f,k}})+\norm{\mc B L_T^{-1} \mc B^{\mathsf{T}}}+\max_{j \in \N_{[1,N_c]}} K_{i,v_j}$, $\tilde{\gamma}_4=\norm{(Y_f-K_{p,v}) C_f^{-1}}$.
\end{condition}
Before discussing the practical implications of Condition \ref{cond.stab} - \ref{cond:current.loop.gains} we present the following theorem that establishes almost global asymptotic stability of the multi-converter system.
\begin{theorem}{\bf(Almost global stability of $\boldsymbol{\mc X^s(\mc S \cap \mc A)}$)}\label{thm.main}
	Consider set-points $p^\star_k$, $q^\star_k$, $v^\star_k$, steady-state angles $\theta^\star_{jk}$, a stability margin $c_1, \ c_2$ and $c_3$, and control gains $\alpha$, $\eta$, $K_{p,i_{\! f,k}}$, $K_{i,i_{\! f,k}}$, $K_{p,v_k}$ and $K_{i,v_k}$ such that Conditions~\ref{cond.stab} - \ref{cond:current.loop.gains} hold. Then, the dynamics \eqref{eq:closed.loop.r.f.dynamcics} are almost globally asymptotically stable with respect to $\mc X^s(\mc S \cap \mc A)$, and the origin $\mathbbl{0}_{n}$ is an exponentially unstable equilibrium.\end{theorem}
The proof of Theorem \ref{thm.main} is given in Appendix~\ref{app:proofs} and is obtained by applying Theorem \ref{theorem:SPAGAS} and Theorem \ref{theorem:SPROA} to the nested dynamics \eqref{eq:closed.loop.r.f.dynamcics} and the Lyapunov functions given in Appendix~\ref{app:lyapred}.
Finally, the next result provides the basis for an applied interpretation of Theorem \ref{thm.main} and Conditions \ref{cond.consistent} - \ref{cond:current.loop.gains}.
\begin{corollary} {\bf (Stabilizing control gains {and set-points})}\label{coro:existence.of.the.gains}
	If the graph $\mc G$ is connected there exist set-points $p_k^{\star}$, $q_k^{\star}$, $v_k^{\star}$, { corresponding} relative nominal phase angles $\theta^{\star}_{jk}$, control gains $\eta_a$ and $\eta$, stability margins $c_1$, $c_2$, $c_3$, and control gains $K_{p,v}$, $K_{i,v}$, $K_{p,f}$, and $K_{i,f}$ such that the dynamics \eqref{eq:closed.loop.r.f.dynamcics} are almost globally asymptotically stable with respect to $\mc X^s(\mc S \cap \mc A)$.
\end{corollary}
The proof is provided in Appendix~\ref{app:proofs}. 

\subsection{Interpretation of the stability conditions}
\begin{figure}
 \includegraphics[width=0.98\columnwidth]{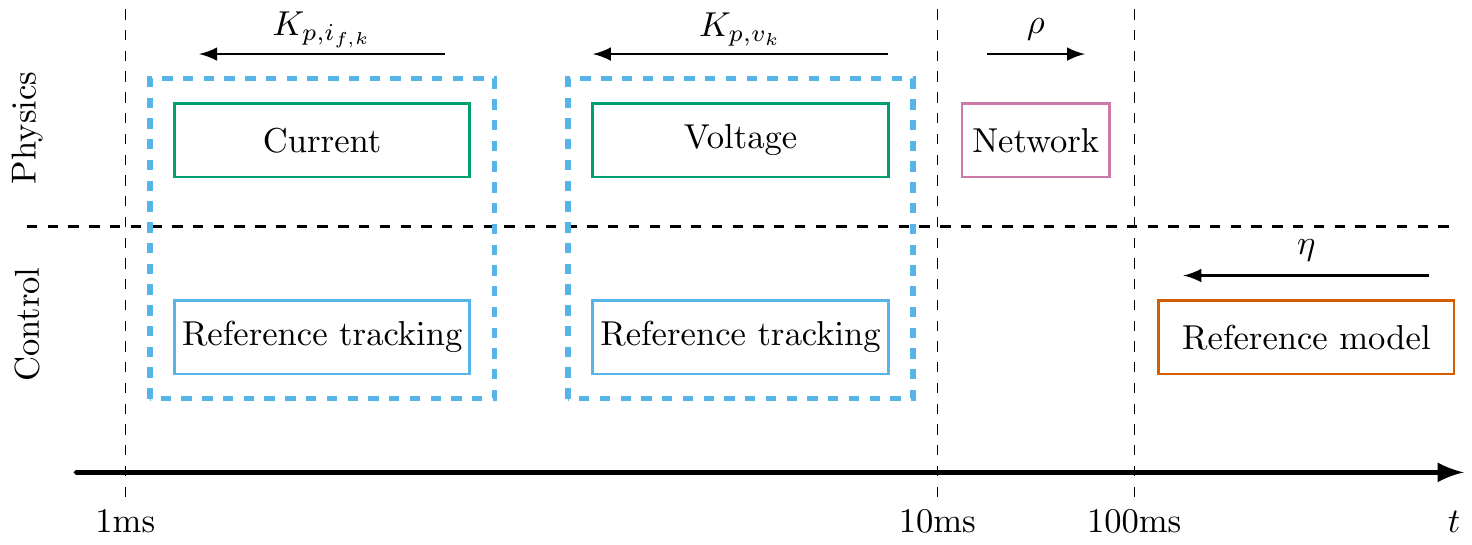}
 \caption{Interpretation of Condition \ref{cond.stab} - \ref{cond:current.loop.gains}. The reference model needs to be sufficiently slow, while the controlled filter current and voltage need to be sufficiently fast. \label{fig:Interpretation}}
\end{figure}

An interpretation of the stability conditions is shown in Figure \ref{fig:Interpretation}. The parameter $\rho$ is the time constant of the transmission lines that cannot be influenced through control. Moreover, the synchronization gain $\eta$ is the time constant of the reference model, and the control gains $K_{p,f}$ and $K_{p,v}$ dominantly influence the convergence rate of the filter current and the terminal voltage closed loop dynamics. Hence, the control gains can be selected such that a sufficient time-scale separation is be enforced and stability is achieved. 

Broadly speaking, the second inequality of Condition \ref{cond.stab} requires the reference model to be slow enough compared to the line dynamics, Condition \ref{cond:voltage.loop.gains} implies that the controlled voltages settle sufficiently fast compared to the transmission line dynamics, and Condition \ref{cond:current.loop.gains} requires the controlled filter current to converge faster than the terminal voltages. 

Specifically, a larger proportional control gain (i.e., speeding up the controlled filter dynamics) and a smaller integral gain (i.e., inducing less oscillations) can be used to 
enforce Condition \ref{cond:voltage.loop.gains}. In particular, a higher proportional gain is required for lower stability margins $c_2$. To this end, we note that
\begin{align*}
0 < c_2 \leq \frac{c_L \rho \norm{\mc \Yn}}{5 \eta \norm{\mc K-\mc L}}
\end{align*}
i.e., the stability margin $c_2$ can be expected to decrease with increased network loading (i.e., smaller $c_L$ and $\norm{\mc K} \to \infty$), faster reference dynamics (i.e., larger $\eta$), and decreased network admittance $\norm{\mc \Yn}$ or time constants $\rho$. Using $\norm{\mc \Yn}=\cos(\kappa)\norm{\mc B R^{-1}_T \mc B^\mathsf{T}} $, the stability margin $c_3$ can be bounded by
\begin{align*}
  0 < c_3 \leq& \frac{1}{2} + \frac{ c_L \sin(\kappa)\beta_{3,1} }{10 \, \omega_0  \norm{\mc K-\mc L} }
\end{align*}
and is therefore expected to decrease with increased network loading (i.e., smaller $c_L$ and $\norm{\mc K} \to \infty$), decreased network inductance (i.e., $\kappa \to 0$), and increased voltage loop control gains (i.e., decreased $\beta_{3,1}$). In turn, $\beta_{3,1}$ is expected to decrease for faster reference dynamics (i.e., larger $\eta$) to enforce Condition \ref{cond:voltage.loop.gains}.

All of these scenarios potentially require larger proportional control gain for the current controller to enforce \ref{cond:current.loop.gains} and a sufficient time scale separation. This reflects practical engineering insights. In particular, the usual power electronics control design begins with tuning the inner current loop and proceeds to the outer loops while enforcing a significant time-scale separation. 

Finally, we note that the conditions \ref{cond.stab} to \ref{cond:current.loop.gains} reflect the nested structure of the system. In particular, given the network parameters, control gains and set-points for the reference model can be computed that satisfy Condition \ref{cond.stab}. Furthermore, given the network parameters, the control gains and set-points of the reference model, and the filter capacitance, stabilizing voltage control gains can be found using Condition \ref{cond:voltage.loop.gains}. Finally, for fixed voltage control gains, the current control gains can be found using Condition \ref{cond:current.loop.gains}.

\begin{remark}[Conservativeness and practical relevance]
For the particular example provided in this section, the numerical bounds provided by Conditions \ref{cond.consistent}-\ref{cond:current.loop.gains} are conservative ($\approx$ a factor of 10 for each timescale). Corollary~\ref{coro:existence.of.the.gains} and the stability conditions  \ref{cond.consistent} - \ref{cond:current.loop.gains} are not intended to provide tight stability ranges but rather as a valuable tuning guideline for the control gains $K_{p,v}$, $K_{i,v}$, $K_{p,f}$, and $K_{i,f}$ by directing the designer on the relative order and direction of tuning, with the guarantee that a stabilizing controller with the given architecture for the full nonlinear system is indeed achievable. 
\end{remark}

\section{Illustrative Example} \label{sec:illustrative.example}
\begin{figure}[t!!!]
	\centering
	\includegraphics[width=.8\columnwidth]{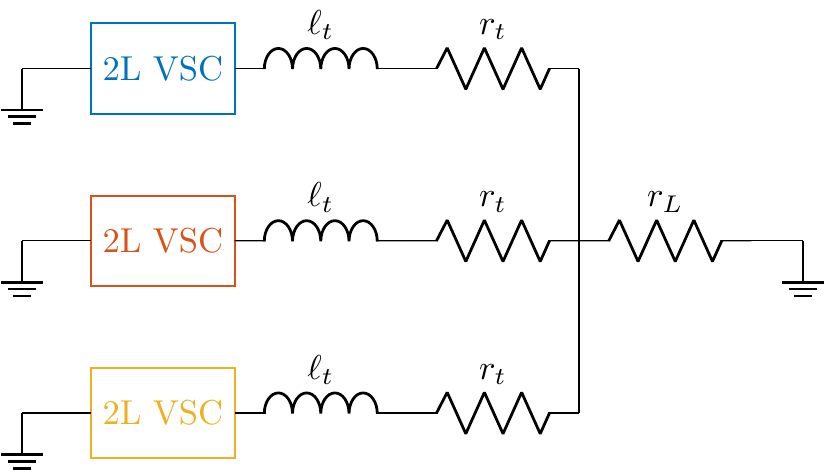}
	\vspace{1em}
	\caption{Three two-level voltage source converters (see Figure \ref{fig:2LVSC}) connected to a resistive load. \label{fig:inverters}}
	\vspace{1em}
	\includegraphics[width=.98\columnwidth]{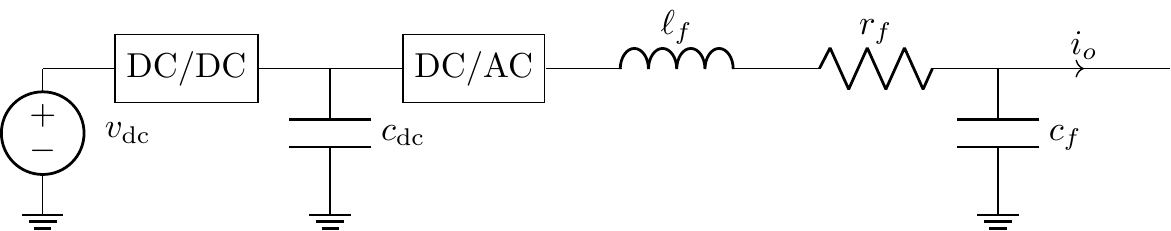}
	\vspace{1em}	
	\caption{Two-level DC/AC voltage source converter with $RLC$ output filter, DC-link capacitor, and a DC/DC boost converter used as power supply and to stabilize the DC-link voltage. \label{fig:2LVSC}}
\end{figure}

As an illustrative example, we use a high-fidelity PLECS \cite{PLECS} model of the single-phase\footnote{The $\beta$-components of the three-phase signals are reconstructed using a Hilbert transform \cite{GC+19}.} microgrid hardware testbed described in \cite{GC+19} with a resistive load shown in Figure \ref{fig:inverters}. In contrast to the theoretical analysis, the switching stage of the voltage source converters in the PLECS model are not averaged and the DC voltage is not constant. Specifically, the simulation uses the model shown in Figure \ref{fig:2LVSC} (see Table \ref{tab:simulation.parameters} for the parameters) that consists of a two-level DC/AC voltage source converter, an RLC output filter, a DC-link capacitor, and a DC/DC boost converter used to stabilize the DC voltage. Moreover, the converter switches are driven via pulse width modulation (with $30$ kHz base frequency). Finally, the controllers implemented in the stationary $(\alpha \beta)$ coordinates \cite{line_dynamics}, and discretized at a frequency of $15$ kHz, and the electromagnetic dynamics are simulated using a variable step ODE solver with a maximum step size of $1.66$ $\mu$s (i.e., twenty times the PWM base frequency). The base power is $S_b=1$ kW, the base voltage is $V_b=120$ V, and, in per unit, the dVOC gains $\eta$ and $\eta_\alpha$ correspond to an active power to frequency droop gain $m_p=\eta S_b/(\omega_0 V^2_b) = 0.15\%$ and reactive power to voltage droop gain $m_q=S_b/(\eta_a V^2_b)=2.2\%$ \cite{GC+19}. While $m_p=0.2\%$ may seem small, it is a realistic choice for droop control of microgrids with low converter power rating \cite[Fig. 4]{vorobev2017high} and short cables \cite[Fig. 3]{vorobev2017high}. Specifically, for low power ratings $S_b$ a larger change in frequency (i.e., $m_p \in [2\%,5\%]$ per unit) for a unit change in power would require a larger gain $\eta$ that can destabilize the system (see Condition \ref{cond.stab}).

\begin{figure*}[t!!!]
	\centering
	\includegraphics[width=0.99\linewidth]{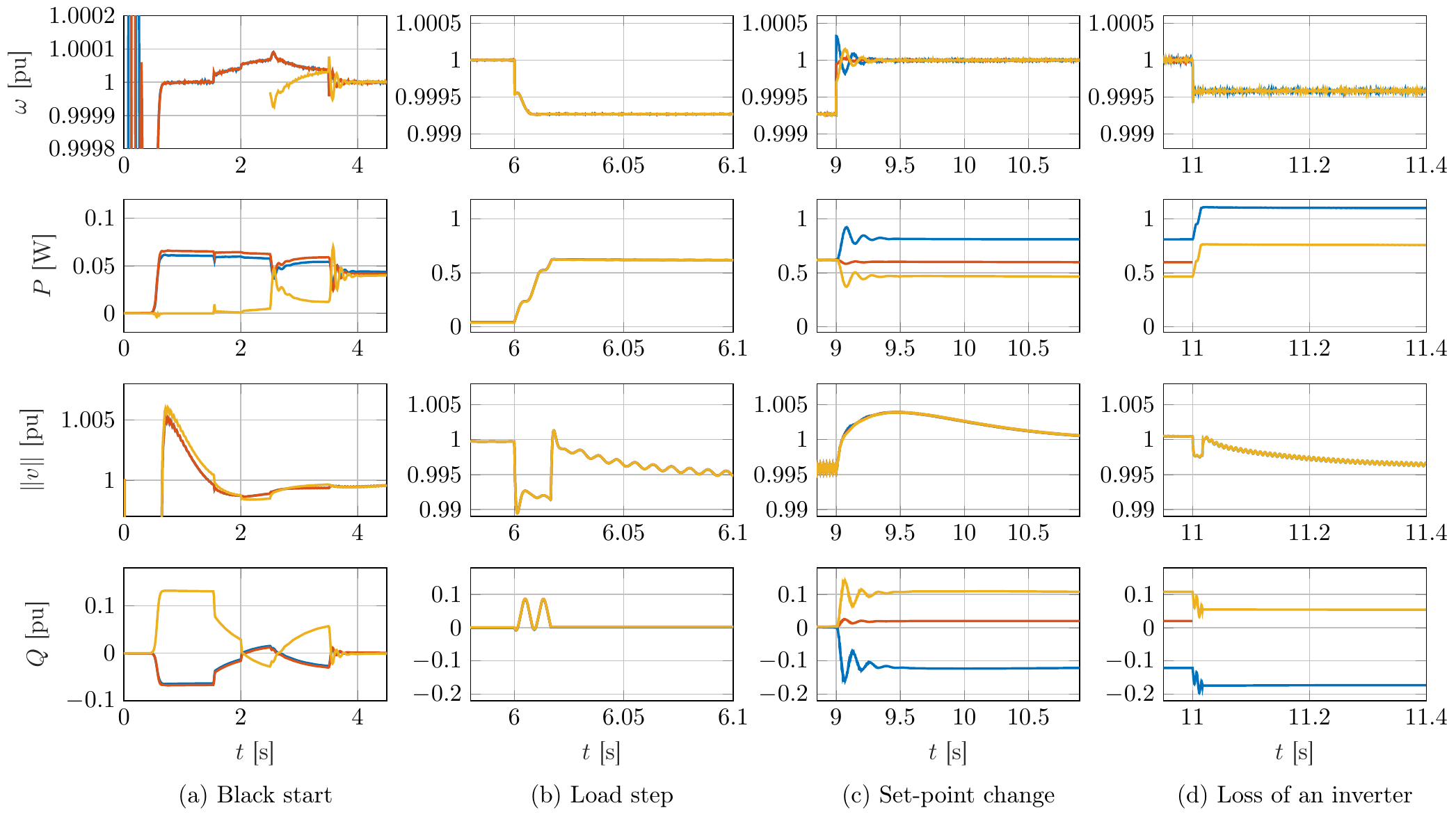}
	\caption{Simulation results showing the frequency, active power, terminal voltage and reactive power for different events. The colors correspond to the converters shown in Figure \ref{fig:inverters}.}
	\label{fig:sim}
\end{figure*}

In order to illustrate the behavior of a closed-loop system for different scenarios, we simulate the sequence of events, shown in Figure \ref{fig:sim}.
\begin{table}[t!]
\begin{center}
\caption{Simulation parameters \label{tab:simulation.parameters}}
\begin{tabular}{lr}
\toprule 
2L-VSC & $r_f=0.124\Omega$, $\ell_f=1$mH, $c_f=24\mu$F, $c_\text{dc}$=1.35mF\\[.3em]
Network & $\omega_0=60$Hz, $r_t=50$m$\Omega$, $\ell_t=0.2$mH\\[.3em]
Load & $r_0$=115$\Omega$,  $r_1$=57.6$\Omega$, $r_2=28.8\Omega$, $r_3=14.4\Omega$ \\
\midrule
  dVOC & $\eta=8.14\tfrac{\text{rad}}{\text{s}}\tfrac{\text{V}}{\text{A}}$, $\eta_a=3.13\tfrac{\text{A}}{\text{V}}$ \\[.3em]
Voltage PI& $K_{p,v}=0.07\tfrac{\text{A}}{\text{V}}$, $K_{i,v}=0.15\tfrac{\text{rad}}{\text{s}}\tfrac{\text{A}}{\text{V}}$ \\[0.3em]
Current PI & $K_{p,f}= 5.93\tfrac{\text{V}}{\text{A}}$, $K_{i,f}=12.49\tfrac{\text{rad}}{\text{s}}\tfrac{\text{V}}{\text{A}}$ \\
  \bottomrule
\end{tabular}

\end{center}
\end{table}
Figure \ref{fig:sim} (a) shows the black start of the microgrid using the first and second converter (for $v^\star_1=v^\star_2=120$V, $p_1^\star= 60.35$W, $q_1^\star= -63.6$var, $p_2^\star= 65$W, and $q_2^\star= -66.7$var). In this setup, the RLC filter of the converters cannot be disconnected, i.e., while the third converter is not operating its output filter acts as reactive power load. Next, between $t=1.5s$ and $t=3.5s$ the third converter is pre-synchronized. First, at $t=1.5$s the third converter is activated using only the current PI controller \eqref{eq:filter.current.controller.single} with reference $i^r_{f,k}=Y_{f,k} v_k$ (i.e., as a current source supplying the filter losses). Next, at $t=2$s the reference model \eqref{eq.closed.loop.r.f.hatv} and voltage PI controller \eqref{eq:voltage.controller.single} are enabled (with $v^\star_3=120V$ and $p^\star_3=q^\star_3=0$), and the feed-forward term $\underline{i}_{o,k}$ in the voltage controller is enabled at $t=2.5$s. Finally, at time $t=3.5$s the set points are changed such that all converters provide equal power to the load, and the frequency, terminal voltage, active, and reactive powers all reach their set-points (i.e., $p_1^\star= 43.2$W, $q_1^\star = -0.97$var, $p_2^\star= 41$W, $q_2^\star= -0.5$var, $p_3^\star=41$V, and $q_3^\star=0.5$var).

A load increase from $125$ W to $1875$ W (without updating the set points) is shown in Figure \ref{fig:sim} (b) illustrating power sharing and the droop like behavior of dVOC.  Figure \ref{fig:sim} (c) shows the response of the converters to a set-point update to $p_1^\star=809.44$W, $q^\star_1=-120.038$var, $p^\star_2=597.12$W, $q^\star_2=19.58$var, $p^\star_3=464.53$W, and $q^\star_3=107.019$var. Finally, Figure \ref{fig:sim} (d) shows the loss of the second converter. Again, the power sharing and droop like behavior of the frequency and active power can be observed. Moreover, the first and the third converter again provide reactive power to the output filter of the second converter.

\section{Conclusion and outlook}\label{sec.conclusion}

In this work, we developed a Lyapunov function framework for stability
analysis of nested nonlinear dynamical systems on multiple time scales where each of the dynamics are influenced only by the slower and the next fastest dynamics. Exploiting this structure and under several technical assumptions (i.e., the existence of steady-state maps and Lyapunov candidate functions used to bound the convergence rates and coupling between the dynamics) we obtain conditions for (almost) global asymptotic stability with respect to a set. Our approach explicitly considers multiple time scales, convergence rates instead of scalar time constants, and reduces conservatism.

As an illustration of our technical result, we apply it to a multiple-converter power system model that includes transmission lines dynamics, with constant inductance to resistance ratio, converter dynamics, cascaded current and voltage control loops, and dVOC as a reference model,  Finally, we obtain explicit stability conditions on the control gains that enforce the well-known time scale separation between the different dynamics, i.e., the dVOC reference model has to be sufficiently slow relative to the line dynamics and the controlled converter voltage and current have to be sufficiently fast compared to the line dynamics. Moreover, the converter current has to be faster than the terminal voltage. Finally, we used a high-fidelity simulation with detailed converter models (i.e., full switching and DC side dynamics) to validate the performance of the proposed control strategy.

This work opens numerous further research directions. For instance, robustness analysis in the presence of parameteric uncertainty are of interest both power systems application to, e.g., account for operating points that are not consistent with the power flow equations and inhomogeneous inductance to resistance ratios, and in the general Lyapunov framework. Also, more general interconnections between the dynamics with different convergence rates should be investigated.

\section*{Acknowledgments}
The authors would like to thank Brian Johnson and Gab-Su Seo for providing the simulation model used in Section \ref{sec:illustrative.example}.

 \begin{appendices}
 
 \section{Lyapunov functions for the multi-converter system}~\label{app:lyapred}
In this appendix, we present the individual Lyapunov function candidates and comparison functions for the reduced-order system $\ddt \hat{v}=f^s_{\hat{v}}(\hat{v})$, the transmission lines, terminal voltages, and filter currents of the multi-converter system discussed in Section~\ref{sec:model} and Section~\ref{sec:control_design}.

Given the voltage set-points $v^\star_k$ and relative steady-state angles $\theta^\star_{k1}$ for all $k \in \mathbb{N}_{[1,N]}$, we define the matrix $S\coloneqq[v_1^\star R(\theta_{11}^\star)^\mathsf{T} \hdots v_N^\star R(\theta_{1N}^\star)^\mathsf{T}]^\mathsf{T}$ whose null space encodes $v_k/v_k^{\star} = R(\theta_{k1}^{\star}) v_1/v_1^{\star}$ for all $k \in \N_{[2,N_c]}$, and the projector $P_S \coloneqq I_{2N} - \frac{1}{\sum v_i^{\star 2}}S S^\mathsf{T}$ onto the nullspace of $\mc S$. Then, the Lyapunov function $V_{\hat{v}}: \R^{2N_c} \to \R_{\geq 0}$ for the reduced-order dynamics \eqref{eq:ss.f.xv} is given by \cite[Prop. 3]{line_dynamics}
\begin{align}\label{eq.v}
V_{\hat{v}}(\hat{v}) \coloneqq \frac{1}{2} \hat{v}^\mathsf{T}P_S\hat{v} + \frac{1}{2}  \eta\eta_a\alpha_1\sum_{k=1}^{N_c} \left( \frac{{v_k^\star}^2 - \norm{\hat{v}_k}^2}{v_k^\star}\right)^2,
\end{align}
where $\eta_a\in \mathbb{R}_{>0}$ is the voltage controller gain and, given $c_L \in \mathbb{R}_{>0}$, the constant $\alpha_1 $ is given by
\begin{align} \label{thm1.strictbound}
\alpha_1 \coloneqq \frac{ c_L  }{ 5\eta \norm{\mc K -\mc L}^2}.
\end{align}
Furthermore, the constants $\eta$ and $c$ cannot be chosen arbitrarily. They must be chosen such that Condition \ref{cond.stab} is always satisfied.
Moreover, we define the comparison function $\psi_{\hat{v}}: \R^{2N_c} \to \R_{\geq 0}$ as $\psi_{\hat{v}}(\hat{v}) \coloneqq \eta  \left(\|\mc K-\mc L\|  \|\hat{v}\|_{\mc S}  + \eta_a \|\Phi(\hat{v})\hat{v}\| \right)$.

Next, let $\mc B_{n} \coloneqq B_n \otimes I_2 \in \R^{N_t \times N_{t_0}}$, where the columns
of the matrix $B_n$ span the nullspace of $B$. Because the graph $\mc G$ is connected it follows from the rank-nullity theorem that $N_{t_0} \coloneqq N_t-N_c+1$. Using $P_t \coloneqq \rho \big( \mc B^{\mathsf{T}} \mc B + L_T \mc B_n \mc B_n^{\mathsf{T}} L_T \big)$ and the error coordinates $y_t \coloneqq i_t - \phi_t(\hat{v})$ we define the candidate Lyapunov function $V_t (y_t)  \coloneqq \frac{1}{2} y^\mathsf{T}_t P_t y_t$ for the transmission lines and the comparison functions $\psi_t(y_t) = \norm{\mc B y_t}$, $\psi^\prime_t(y_t) =  \norm{\mc B_n L_T y_t}$.
\begin{remark}[Circulating network currents]\label{rem:circcurr}
The network currents $i_t$ can separated into current injections $\mc B i_t$ and circulating currents $\mc B_n i_t$ that do not affect the converters. Using $\psi_t=\norm{\mc B y_t}$ and $\psi^\prime_t=\norm{\mc B_n L_t y_t}$ these currents can be separated in the bounds of Assumption \ref{assumption:derivative}, i.e., only the currents $\mc B_n i_t$ enter into the bounds on the coupling between different time scales.
\end{remark}
Moreover, using the error coordinates $y_v\coloneqq x_v-\phi_{v}$ and 
\begin{align*} 
P_v\coloneqq \begin{bmatrix}
 K_{p,v}^{-1} C_f & K_{i,v}^{-1}C_f \\  K_{i,v}^{-1}C_f & K_{p,v}K_{i,v}^{-1}+K_{i,v}K_{p,v}^{-1}
\end{bmatrix}
\end{align*}
we define the candidate Lyapunov function $V_v(y_v) \coloneqq \frac{1}{2} y_v^{\mathsf{T}} P_v y_v$  for the terminal voltages. Finally, using the error coordinates $y_f  \coloneqq  x_f -  \phi_{f}$ and
\begin{align*}
 P_f \coloneqq  \begin{bmatrix}
K_{p,f}^{-1} L_f & K_{i,f}^{-1}L_f \\ K_{i,f}^{-1}L_f & K_{p,f} K_{i,f}^{-1} + K_{i,f}K_{p,f}^{-1}
\end{bmatrix}
\end{align*}
we define the candidate Lyapunov function $V_{f} \coloneqq \frac{1}{2} y_f^{\mathsf{T}} P_f y_f$ for the filter currents.

\section{Proofs and technical lemmas}\label{app:proofs}

{\begin{lemma}{\bf(Lyapunov decrease)}\label{lemma:decrease.epsilon} For all $i \in  \N_{[1,N]}$ let
$\alpha^{(\epsilon)}_i \coloneqq \epsilon_i^{-1} \alpha_i$, for all $i \in  \N_{[1,N-1]}$ let $\beta^{(\epsilon)}_{i,i+1} \coloneqq \epsilon_i^{-1} \beta_{i,i+1}$, for all $i \in  \N_{[2,N]}$ let $\beta^{(\epsilon)}_{i,1} \coloneqq \sum\nolimits_{k=1}^{i-1}\epsilon_k^{-1} b_{i,1,k}$ and $\gamma^{(\epsilon)}_i \coloneqq \epsilon_{i-1}^{-1}b_{i,i,i-1}$,  and for all $i \in  \N_{[3,N]}$ and all $j \in \N_{[2,i-1]}$ let $\beta^{(\epsilon)}_{i,j} \coloneqq \sum\nolimits_{k=j-1}^{i-1} \epsilon_k^{-1} b_{i,j,k}$. Consider the function $\nu_\epsilon = \mu^{(\epsilon)}_1 V_1(x_1) + \sum_{i=2}^{N} \mu^{(\epsilon)}_i V_i (y_i)$ with $\mu^{(\epsilon)}_i = \prod_{j=1}^{i-1} \frac{\beta^{(\epsilon)}_{j,j+1}}{\beta^{(\epsilon)}_{j+1,j}}$ for all $i \in \N_{[1,N]}$, and $\mu^{(\epsilon)}_1=1$. Under Assumption \ref{assumption:derivative}, the derivative of $\nu_\epsilon$
along the trajectories of the \eqref{eq.epsilon.existence}
	is bounded by
\begin{align*}
  \ddt \nu_\epsilon \leq - \left[ \begin{smallmatrix} \psi^{\phantom{\prime}}_1 \\ \psi^{\phantom{\prime}}_2 \\ \vdots \\ \psi^{\phantom{\prime}}_N \end{smallmatrix} \right]^{\mathsf{T}} \!\!H \left[ \begin{smallmatrix} \psi^{\phantom{\prime}}_1 \\ \psi^{\phantom{\prime}}_2\\ \vdots \\ \psi^{\phantom{\prime}}_N  \end{smallmatrix} \right] - \left[ \begin{smallmatrix} \psi^\prime_1 \\ \psi^\prime_2 \\ \vdots \\ \psi^\prime_N \end{smallmatrix} \right]^{\mathsf{T}} \!\!H^\prime \left[ \begin{smallmatrix} \psi^\prime_1 \\ \psi^\prime_2\\ \vdots \\ \psi^\prime_N  \end{smallmatrix} \right],
\end{align*}
where $H^\prime\coloneqq\diag(\{\mu^{(\epsilon)}_i \alpha^{(\epsilon)\prime}_i\}_{i=1}^{N})$, and $H$ is defined recursively for all $i \in \N_{[2,N]}$ starting from $H_1=\alpha^{(\epsilon)}_1$ by
\begin{align*}
H_i=
\left[
\begin{array}{c|c}
H_{i-1} & \boldsymbol{\beta}^{(\epsilon)}_i \mu^{(\epsilon)}_i \\
\hline
\star
& (\alpha^{(\epsilon)}_i-\gamma^{(\epsilon)}_i)\mu^{(\epsilon)}_i
\end{array} \right],
\end{align*}
and $\boldsymbol{\beta}^{(\epsilon)}_i \coloneqq (\ldots,
\frac{1}{2}\beta^{(\epsilon)}_{i,i-3},\frac{1}{2}\beta^{(\epsilon)}_{i,i-2},\beta^{(\epsilon)}_{i,i-1})$.
\end{lemma}
\begin{IEEEproof}
From the time derivative of \eqref{eq.muLyap} we conclude that
\begin{subequations}
\begin{align*}
\ddt{\nu_\epsilon} \coloneqq&  \frac{\mu^{(\epsilon)}_1}{\epsilon_1} \dxdy{V_1}{x_1} f_1 
+ \sum\nolimits_{i=2}^{N-1} \frac{\mu^{(\epsilon)}_i}{\epsilon_i} \dxdy{V_i}{y_i} {f_i}
+\frac{\mu^{(\epsilon)}_N}{\epsilon_N} \dxdy{V_N}{y_N} f_N \\
&- \sum\nolimits_{i=2}^{N} \sum\nolimits_{k=1}^{i-1} \frac{\mu^{(\epsilon)}_i}{\epsilon_k} \dxdy{V_i}{y_i} \dxdy{\phi_i}{x_k} f_k 
\end{align*}
\end{subequations}
Next, we add and subtract 
$\frac{\mu^{(\epsilon)}_i}{\epsilon_i}\dxdy{V_i}{y_i}  f^s_i$ for all $i \in \N_{[1,N-1]}$.
Using Assumption \ref{assumption:derivative} one obtains
\begin{align*}
	\frac{\mu^{(\epsilon)}_1}{\epsilon_1} \dxdy{V_1}{x_1}f^s_1 
	&\!\leq - \mu^{(\epsilon)}_1 \alpha^{(\epsilon)}_1 \psi_1(x_1)^2 - \mu^{(\epsilon)}_1 \alpha^{(\epsilon)\prime}_1 \psi^{\prime}_1(x_1)^2,\\
	\frac{\mu^{(\epsilon)}_1}{\epsilon_1}  \dxdy{V_1}{x_1}\big(f_1 
	\!-\! f^s_1 
	\big) &\!\leq  \mu^{(\epsilon)}_1 \beta^{(\epsilon)}_{12}\psi_1(x_1)\psi_2(y_2),\\
	\frac{\mu^{(\epsilon)}_N}{\epsilon_N} \dxdy{V_N}{y_N} f_N 
	&\!\leq \!\! -\mu^{(\epsilon)}_N \alpha^{(\epsilon)}_N \psi_N(y_N)^2 -\mu^{(\epsilon)}_N \alpha^{(\epsilon)\prime}_N {\psi_N^{\prime}(y_N)}^2\!,
\end{align*}
and for all $i \in \N_{[2,N-1]}$ one obtains
\begin{align*}
	  \frac{\mu^{(\epsilon)}_i}{\epsilon_i} \dxdy{V_i}{y_i}{f^s_i}
	  &\leq - \mu^{(\epsilon)}_i \alpha^{(\epsilon)}_i \psi_i(y_i)^2 - \mu^{(\epsilon)}_i \alpha^{(\epsilon)\prime}_i {\psi_i^{\prime}(y_i)}^2,\\
	  \frac{\mu^{(\epsilon)}_i}{\epsilon_i}\dxdy{V_i}{y_i}\bigg (f_i
	  -f^s_i 
	  \bigg)  &\leq \mu^{(\epsilon)}_i \beta^{(\epsilon)}_{i,i+1} \psi_i(y_i)\psi_{i+1}(y_{i+1}).
\end{align*}
	Next, using Assumption \ref{assumption:derivative}, for all $i \in \N_{[2,N]}$ it holds that
	\begin{align*}
	 &-\sum_{k=1}^{i-1} \frac{\mu^{(\epsilon)}_i}{\epsilon_k} \dxdy{V_i}{y_i} \dxdy{\phi_i}{x_k} f_k 
	 \leq \sum_{k=1}^{i-1} \frac{\mu^{(\epsilon)}_i}{\epsilon_k}  \sum_{j=1}^{k+1} b_{i,j,k} \psi_i \psi_j=\\
	 &=\mu^{(\epsilon)}_i \psi_i \bigg( \underbrace{\sum_{k=1}^{i-1} \frac{b_{i,1,k}}{\epsilon_k}}_{=\beta^{(\epsilon)}_{i,1}} \psi_{1} +\sum_{j=2}^{i-1} \underbrace{\sum_{k=j-1}^{i-1}\frac{b_{i,j,k}}{\epsilon_k}}_{=\beta^{(\epsilon)}_{i,j}} \psi_j +\underbrace{\vphantom{\sum_{k=j-1}^{i-1}} \frac{b_{i,i,i-1}}{\epsilon_{i-1}}}_{=\gamma^{(\epsilon)}_i} \psi_i \bigg),
	\end{align*}
	where $\psi_1 = \psi_1(x_1)$ and $\psi_i = \psi_i(y_i)$ for all $i \in \N_{[2,N]}$.
	Using these bounds, we can bound $\ddt \nu_\epsilon$ by two quadratic forms
	\begin{align*}
  \ddt \nu_\epsilon \leq - \left[ \begin{smallmatrix} \psi^{\phantom{\prime}}_1 \\ \psi^{\phantom{\prime}}_2 \\ \vdots \\ \psi^{\phantom{\prime}}_N \end{smallmatrix} \right]^{\mathsf{T}} \!\!H^{(\mu)} \left[ \begin{smallmatrix} \psi^{\phantom{\prime}}_1 \\ \psi^{\phantom{\prime}}_2\\ \vdots \\ \psi^{\phantom{\prime}}_N  \end{smallmatrix} \right] - \left[ \begin{smallmatrix} \psi^\prime_1 \\ \psi^\prime_2 \\ \vdots \\ \psi^\prime_N \end{smallmatrix} \right]^{\mathsf{T}} \!\!H^\prime \left[ \begin{smallmatrix} \psi^\prime_1 \\ \psi^\prime_2\\ \vdots \\ \psi^\prime_N  \end{smallmatrix} \right],
\end{align*}
where $H^\prime\coloneqq\diag(\{\mu^{(\epsilon)}_i \alpha^{(\epsilon)\prime}_i\}_{i=1}^{N})$ and $H$ is defined recursively for all $i \in \N_{[2,N]}$ starting from $H^{(\mu)}_1=\alpha^{(\epsilon)}_1$ by
\begin{align*}
H^{(\mu)}_i&=%\scalebox{0.8}{$ 
\left[
\begin{array}{c|c}
\multirow{4}{*}{$H^{(\mu)}_{i-1}$} &\vdots \\[0.2em]
& -\tfrac{1}{2} \mu^{(\epsilon)}_i \beta^{(\epsilon)}_{i,i-3}  \\[0.2em]
& -\tfrac{1}{2} \mu^{(\epsilon)}_i \beta^{(\epsilon)}_{i,i-2}  \\[0.2em]
& -\tfrac{1}{2} \mu^{(\epsilon)}_i \big(\tfrac{\mu^{(\epsilon)}_{i-1}}{\mu^{(\epsilon)}_i} \beta^{(\epsilon)}_{i-1,i} + \beta^{(\epsilon)}_{i,i-1}\big)  \\[0.4em]
\hline
\star
&  (\alpha^{(\epsilon)}_i-\gamma^{(\epsilon)}_i) \mu^{(\epsilon)}_i
\end{array} \right].
%$}
\end{align*}
The lemma follows by noting that $\tfrac{\mu^{(\epsilon)}_{i-1}}{\mu^{(\epsilon)}_i} \beta^{(\epsilon)}_{i-1,i} = \beta^{(\epsilon)}_{i,i-1}$.	
\end{IEEEproof}}
	
	\begin{lemma}{\bf(Positive definiteness)} \label{lemma:positive.definiteness}
		Consider a symmetric matrix $A=\{a_{i,j} \}_{N \times N}$, with $a_{i,i}>0$ for all $i \in \mathbb{N}_{[1,N-1]}$ and $a_{i,j}=a_{j,i} \leq 0,$ for all $i \in \N_{[1,N]}$,  $j \in \N_{[1,N]}$, and $i \neq j$. For all $k \in \mathbb{N}_{[1,N]}$, let $A_k=\{a_{i,j} \}_{k \times k}$ denote the $k$-th leading principal minor of $A$. $A$ is positive definite if and only if for all $k \in \mathbb{N}_{[2,N]}$, $A_k$ is invertible and satisfies
%		%
		\begin{align}
		a_{k,k}>\left[\begin{smallmatrix}
		a_{k,1} \\ \vdots \\ a_{k,k-1} 
		\end{smallmatrix} \right]^{\mathsf{T}} A_{k-1}^{-1} \left[\begin{smallmatrix}
		a_{1,k} \\ \vdots \\ a_{k-1,k} 
		\end{smallmatrix} \right] \geq 0 \label{eq:lemma}, 
		\end{align}
		where equality holds if $\left[\begin{smallmatrix}
		a_{k,1} \\ \vdots \\ a_{k,k-1} 
		\end{smallmatrix} \right]=\left[\begin{smallmatrix}
		a_{1,k} \\ \vdots \\ a_{k-1,k} 
		\end{smallmatrix} \right]= \mathbbl{0}_{k-1}$.
	\end{lemma}
	\begin{IEEEproof}
		Using the Schur complement, it can be verified that the $k$-th  minor of $A$ is positive definite if and only if 
		\begin{align*}
		a_{k,k}-\left[\begin{smallmatrix}
		a_{k,1} \\ \vdots \\ a_{k,k-1} 
		\end{smallmatrix} \right]^{\mathsf{T}} A_{k-1}^{-1} \left[\begin{smallmatrix}
		a_{1,k} \\ \vdots \\ a_{k-1,k} 
		\end{smallmatrix} \right] > 0
		\end{align*}
		and $A_{k-1}$ is positive definite. In particular, $A_{N-i}$ is positive definite if and only if \eqref{eq:lemma} holds for $k=N-i$ and $A_{N-1-i}$ is positive definite. The Lemma follows by induction over $i \in \N_{[0,N-2]}$.
	\end{IEEEproof}\label{key}
{\emph{Proof of Theorem \ref{theorem:epsilon}:}
Consider the Lyapunov function candidate $\nu_\epsilon = \mu^{(\epsilon)}_1 V_1(x_1) + \sum_{i=2}^{N} \mu^{(\epsilon)}_i V_i (y_i)$ with $\mu^{(\epsilon)}_i = \prod_{j=1}^{i-1} \frac{\beta^{(\epsilon)}_{j,j+1}}{\beta^{(\epsilon)}_{j+1,j}}$ for all $i \in \N_{[1,N]}$, and $\mu^{(\epsilon)}_1=1$. Using Lemma \ref{lemma:decrease.epsilon} we obtain
\begin{align*}
  \ddt \nu_\epsilon \leq - \left[ \begin{smallmatrix} \psi^{\phantom{\prime}}_1 \\ \psi^{\phantom{\prime}}_2 \\ \vdots \\ \psi^{\phantom{\prime}}_N \end{smallmatrix} \right]^{\mathsf{T}} \!\!H \left[ \begin{smallmatrix} \psi^{\phantom{\prime}}_1 \\ \psi^{\phantom{\prime}}_2\\ \vdots \\ \psi^{\phantom{\prime}}_N  \end{smallmatrix} \right] - \left[ \begin{smallmatrix} \psi^\prime_1 \\ \psi^\prime_2 \\ \vdots \\ \psi^\prime_N \end{smallmatrix} \right]^{\mathsf{T}} \!\!H^\prime \left[ \begin{smallmatrix} \psi^\prime_1 \\ \psi^\prime_2\\ \vdots \\ \psi^\prime_N  \end{smallmatrix} \right].
\end{align*}
Next, using Lemma \ref{lemma:positive.definiteness} it follows that $H$ is positive definite if 
\begin{align}\label{eq:stab.cond.eps}
\gamma^{(\epsilon)}_i+ \mu^{(\epsilon)}_i {\boldsymbol{\beta}^{(\epsilon)}_i}^{\mathsf{T}} H_{i-1}^{-1} {\boldsymbol{\beta}^{(\epsilon)}_i}^{\mathsf{T}}  < \frac{\alpha_i}{\epsilon_i},
\end{align}
holds for all $i\in \mathbb{N}_{[2,N]}$. Because the left hand side of \eqref{eq:stab.cond.eps} only depends on $\epsilon_j$ for all $j \!<\! i$ it directly follows that, for every $i\!\in\! \mathbb{N}_{[2,N]}$ there exists $\epsilon_i \!\in\! \R_{>0}$ with $\epsilon_i \!<\! \epsilon_{i-1}$ such that \eqref{eq:stab.cond.eps} holds. Using the same steps as in the proof of Theorem \ref{theorem:SPAGAS} there exists a function $\chi_3 \!\in\! \mathscr{K}$ such that $\ddt \nu \leq -\chi_3(\norm{x}_{\mc X^s(\mc C_1)})$ and $\chi_1 \!\in\! \mathscr{K}_{\infty}$ and $\chi_2 \!\in\! \mathscr{K}_{\infty}$ such that $\chi_1(\norm{x}_{\mc X^s(\mc C_1)}) \leq \nu(x) \leq \chi_2(\norm{x}_{\mc X^s(\mc C_1)})$. Letting $\mc V \!=\! \nu_\epsilon$, $\mc C = \mc X^s(\mc C_1)$, and $\mc U = \emptyset$, it follows from Theorem \ref{theorem:AGAS} that the system is almost globally asymptotically stable with respect to $\mc C = \mc X^s(\mc C_1)$. Since $\mc U = \emptyset$, \eqref{eq.AGA} holds for all $x \in \R^n$, and it follows that \eqref{eq.epsilon.existence} is globally asymptotically stable.\hfill\IEEEQED}\\[.5em]
\emph{Proof of Corollary \ref{coro:existence.of.the.gains}:}
It can be verified that Condition \ref{cond.consistent} holds for the operating point $v^\star_k = v^\star_j$ and $\theta^\star_{jk}=0$ for all $(j,k) \in \N_{[1,N_c]} \times \N_{[1,N_c]}$ (i.e., the DC power flow) and Condition 4 reduces to $\eta_a \leq \frac{1}{2} \lambda_2(L) - c_L$ and $\eta  \leq \frac{c_L}{ \rho \norm{\mc \Yn} (c_L  + 5 \| \mc L\|)}$. If the graph $\mc G$ is connected it follows that $\lambda_2(L) \in \R_{>0}$, the first inequality holds for $\eta_a = c_L =\frac{1}{4} \lambda_2(L)$, and the right hand side of the second inequality is strictly positive, i.e., there exists $\eta \in \R_{>0}$ such that Condition \ref{cond.stab} holds. Next, note that the right-hand side of the inequality in Condition \ref{cond:voltage.loop.gains} is positive and independent of the gains $K_{p,v_k}$ and $K_{i,v_k}$. Moreover, letting $\xi_p = K_{i,v_k} / K_{p,v_k} \in \R_{>0}$ and $\xi_i = K_{i,v_k} / c_{f,k} \in \R_{>1}$ for all $k \in \N_{[1,N]}$, the left-hand side of the inequality in Condition \ref{cond:voltage.loop.gains} becomes $\frac{1 + \xi_p}{\xi_i - 1} \in \R_{>0}$ and can be made arbitrarily small by choosing $\xi_p \in \R_{>0}$ small enough and $\xi_i \in \R_{>1}$ large enough. Using the same arguments, and noting $c_\beta$ cancel in the rhs of condition 6, it can be verified that there always exists $K_{p,f_k} \in \R_{>0}$ and $K_{i,f_k} \in \R_{>0}$ such that Condition \ref{cond:current.loop.gains} holds. Thus, Conditions \ref{cond.consistent} - \ref{cond:current.loop.gains} hold and Theorem \ref{thm.main} states the desired result.
 \hfill\IEEEQED\\[.5em]

\emph{Proof of Theorem \ref{thm.main}:
We first apply Theorem \ref{theorem:SPAGAS} to the nested system defined by $N=4$, $(x_1,x_2,x_3,x_4)= (\hat{v},i_t,x_v,x_f)$, and $(f_1,f_2,f_3,f_4)=(f_{\hat{v}},f_t,f_v,f_f)$, with 
 the Lyapunov function candidates $(V_1,V_2,V_3,V_4)=(V_{\hat{v}},V_t,V_v,V_f)$ and the comparison functions $(\psi_1,\psi_2,\psi_3,\psi_4)=(\psi_{\hat{v}},\psi_t,\norm{y_v},\norm{y_f})$, $(\psi^\prime_1,\psi^\prime_2,\psi^\prime_3,\psi^\prime_4)=(0,\psi^\prime_t,0,0)$ defined in Appendix~\ref{app:lyapred}, and the sets $\mathcal{C}_1=\mathcal{S} \cap \mathcal{A}$, and $\mathcal{U}=\{\mathbbl{0}_{2N_c}\}$. Note that Conditions \ref{cond:voltage.loop.gains}-\ref{cond:current.loop.gains} imply that $1-\max_{k\in\N_{[1,N_c]}}  \frac{c_{f,k}}{K_{i,v_k}} > 0$ and $1-\max_{k\in\N_{[1,N_c]}} \frac{\ell_{f,k}}{K_{i,f_k}} > 0$. It follows that the matrices $P_v$ and $P_f$ are positive definite, i.e., there exists $\chi^{\scriptscriptstyle{V_i}}_1 \in \mathscr{K}_{\infty}$ and $\chi^{\scriptscriptstyle{V_i}}_2 \in \mathscr{K}_{\infty}$ such that the conditions of Theorem \ref{theorem:SPAGAS} hold for $i=\{3,4\}$ and, by Condition~\ref{cond.stab} and  \cite[Prop. 3-4]{line_dynamics}, they also hold for $i=\{1,2\}$.} 

Using Condition \ref{cond.stab}, $\alpha_1 \in \R_{>0}$ defined in \eqref{thm1.strictbound}, $\alpha_2=1$, and \cite[Prop. 3-4]{line_dynamics} it can be verified that the decrease condition in Assumption \ref{assumption:derivative} holds for $i=\{1,2\}$. Moreover, using $\alpha_3 = 1-\max_{k \in \N_{[1,N_c]}}\frac{c_{f,k}}{K_{i,v_k}}>0$,
$%\begin{align*}
\dxdy{V_3}{y_3} f_3^s = - y_3^{\mathsf{T}}  \diag( I_{2N_c}-K_{i,v}^{-1}C_f, I_{2N_c}) y_3$, and %\end{align*}
 $\alpha_4 = 1 - \max_{k \in \N_{[1,N_c]}}  \frac{\ell_{f ,k}}{K_{i,f_k}} >0$, 
$%\begin{align*}
 \frac{\partial V_4}{\partial y_4}  f_4^s= - y_4^{\mathsf{T}} \diag(	I_{2N_c}-K_{i,f}^{-1}L_f, I_{2N_c}) y_4$ %\end{align*}
it can be seen that $\dxdy{V_i}{y_i} f_i^s \leq -\alpha_i \psi_i^2$ holds for $i=\{3,4\}$.

Next, we show that for $i \in \N_{[1,3]}$ constants $\beta_{i,i+1}$, that satisfy Assumption \ref{assumption:derivative} exist. For $i=1$ the assumption holds for $\beta_{1,2} \!=\! \norm{\mc K - \mc L}^{-1}$  (see \cite[Prop. 5]{line_dynamics}). For $i=2$, let $\beta_{2,3} = \norm{ \mc B R_T^{-1}\mc B^{\mathsf{T}}}$ and note that $f_2$ is separable and linear in its arguments, it holds that 
%\begin{align*}
$\dxdy{V_{2}}{y_{2}}(f_2\!-\!f_2^s)\!=\!\rho \big( y_2^{\mathsf{T}} \mc B^{\mathsf{T}}\! \mc B L_T^{-1} \!+\! y_2^{\mathsf{T}}\! L_T \mc B_n \mc B_n^{\mathsf{T}} \big ) \mc B^{\mathsf{T}}\! y_3 \!\leq\! \beta_{2,3} \psi_2 \psi_3,$ 
%\end{align*}
where we used $\mc B \mc B_{n}=\mathbbl{0}_{2N_c \times 2N_{t_0}}$. Finally, for $i\!\!=\!\!3$ we obtain $f_3-f^s_3= \diag(C^{-1}_f,\mathbbl{0}_{2N_c \times 2N_c}) y_4$ and $\dxdy{V_3}{y_3}(f_3\!-\!f_3^s) \leq \beta_{3,4} \psi_3 \psi_4$ holds for $\beta_{3,4}  \! \coloneqq \! \max_{k \in \N_{[1,N_c]}}  \tfrac{1}{K_{i,v_k}}+\tfrac{1}{K_{p,v_k}}$.

It remains to show that for $i \in \N_{[2,4]}$, $k \in \N_{[1,i-1]}$ and $j \in \N_{[1,k+1]}$, there exist $b_{i,j,k}$ such that Assumption 2 is satisfied. For $i=2$, we have $b_{2,1,1}\!=\! \rho \norm{\Yn}$, and $b_{2,2,1} \! =\! \eta \rho \norm{\Yn}$ (cf. \cite[Prop. 6]{line_dynamics}). Moreover, for $i=3$, we note that $\dxdy{\phi_3}{x_2}=\mathbbl{0}_{n_3}$, hence, $b_{3,j,2}=0$, for all $j \in \N_{[1,3]}$. Let $b_{3,1,1}= \max_{k \in \N_{[1,N_c]}} \tfrac{c_{f,k}}{K_{i,v_k}}+\tfrac{c_{f,k}}{K_{p,v_k}}$ and $b_{3,2,1} = \eta b_{3,1,1}$.  Next, using $\dxdy{\phi_3}{x_1} f_1=f^s_1(x_1) - \eta \mc R(\kappa) \mc B y_2 $, we have $ - \dxdy{V_3}{y_3}\dxdy{\phi_3}{x_1} f_1\leq b_{3,1,1} \psi_1 \psi_3 + b_{3,2,1} \psi_2 \psi_3$.
where we used that $\norm{(\mc K-\mc L) x_1}=\norm{(\mc K-\mc L) P_S x_1} \leq \norm{\mc K - \mc L} \norm{x_1}_{\mc S}$, i.e., $\norm{f^s_1} \leq \psi_{1}$.

Using separability in all and linearity in some arguments of $f_{1}$, $f_2$ and $f_3$, and $Y_f^\prime=(Y_f-K_{p,v})C_f^{-1}$, it follows that $\dxdy{\phi_{4}}{x_1} f_{1} = \big( K_{p,v} (f_{1}^s(x_1)-\eta \mc R (\kappa) \mc B  y_2),\mathbbl{0}_{2N_c} \big)$, $\dxdy{\phi_{4}}{x_{2}} f_2  \!= \linebreak \! \big (\! \big[ \mc B L_T^{-1} \mc B^{\mathsf{T}} \ \mathbbl{0}_{2N_c \times 2N_c} \big ] y_3 \! - \! \mc B (\tfrac{1}{\rho}I_{n_2} \!\! +\! \omega_0 \mc J_{n_2} )y_2, \mathbbl{0}_{2N_c} \! \big)$, and \linebreak $\dxdy{\phi_{4}}{x_{3}}\!f_{{3}} \!\!=\! \ \big( \! \big[Y_f^\prime \  \mathbbl{0}_{2N_c \times 2N_c}
\big]   y_4 \!\!-\! \big[Y_f^\prime K_{p,v} \! +\!  K_{i,v} Y_f^\prime \ \ K_{i,v}  \big] y_3, \! \mathbbl{0}_{2N_c}  \! \big)$.

Using $c_\beta \!=\! \max_{k \in \N_{[1,N_c]}}\!\!  \tfrac{\ell_{f,k}}{K_{i,f_k}} \!\!+\!\! \tfrac{\ell_{f,k}}{K_{p,f_k}}$, $b_{4,1,1}\!=\!c_\beta \max_{k \in \N_{[1,N_c]}} K_{p,v_k}$, $b_{4,2,1}\!=\!c_\beta \eta \max_{k \in \N_{[1,N_c]}} K_{p,v_k}$, $b_{4,2,2}=\tfrac{c_\beta}{\rho \cos(\kappa)}$, $b_{4,3,2}\!=\!c_\beta\norm{\mc B L_T^{-1} \mc B^{\mathsf{T}}}$, $b_{4,3,3}\!=\!c_\beta\big( \norm{Y_f-K_{p,v}} \max_{k \in \N_{[1,N_c]}}\tfrac{K_{p,v_k}}{c_{f,k}}+\tfrac{K_{i,v_k}}{c_{f,k}} + \max_{k \in \N_{[1,N_c]}} K_{i,v_k} \big)$, $b_{4,4,3}\!=\!c_\beta\norm{(Y_f-K_{p,v})C_f^{-1}}$, $b_{4,1,3}\!\!=\!b_{4,1,2}\!\!=\!b_{4,2,3}\!=\!0$, and $\beta_{3,4}\!=\!\max_{k \in \N_{[1,N_c]}} \!\tfrac{1}{K_{i,v_k}\!}+\tfrac{1}{\!K_{p,v_k}}$ it can be verified that $-\dxdy{V_{4}}{y_{4}}\dxdy{\phi_{4}}{x_{k}}f_{k} \leq  \sum_{j=2}^{k+1}  b_{4,j,k}  \psi_4(y_4) \psi_j(y_j)  +   b_{i,1,k}  \psi_4(y_4) \psi_1(x_1)$ holds for all $k \in \N_{[1,3]}$,

To apply Theorem \ref{theorem:SPAGAS}, note that Conditions \ref{cond.stab} to \ref{cond:current.loop.gains} ensure that Condition 2 holds and it follows from Proposition 1 that the matrix $M$ as in Definition \ref{def.recursive.matrix.M} is positive definite.

Finally, we show that the region of attraction of the equilibrium $x^\star = \mathbbl{0}_n$ has measure zero. To this end, we define $V_{\delta,1} \coloneqq \delta^\mathsf{T}_{x,1} (P_S-2\eta \eta_a\alpha_1 I_{2N}) \delta_{x,1}$. Note that the functions $f^s_2$, $f^s_3$, and $f^s_4$ are linear and the corresponding Lyapunov function candidates $V_2$, $V_3$, $V_4$ are quadratic. Moreover, for all $\delta^\prime_{x,1} \in \mc S \setminus \{\mathbbl{0}_{2N_c}\}$ it holds that $V_{\delta,1}(\delta^\prime_{x,1}) < V_{\delta,1}(\mathbbl{0}_{2 N_c})$. Replacing $\psi_1$ with $\psi_1^{0} \coloneqq \eta(\norm{\mc K - \mc L} \norm{\hat{v}_\delta}_{\mc S} + \eta_a\norm{\hat{v}_\delta})$ and using the same arguments and constants as above it can be verified that Assumption \ref{assumption:derivative} holds for the linearized voltage reference dynamics $\ddt \delta_{x,1} = \eta (\mc K \delta_{x,1} - \mc R (\kappa) \mc B \delta_{x,2} + \eta_a \delta_{x,1})$ and the linear dynamics \eqref{eq:transmission_lines_dynamics_rotating}, \eqref{eq.closed.loop.xv} to \eqref{eq.closed.loop.r.f.xif}. Therefore, the conditions of Theorem \ref{theorem:SPROA} are satisfied and the region of attraction of the origin under \eqref{eq:closed.loop.r.f.dynamcics} has measure zero. Therefore, the conditions of Theorem \ref{theorem:SPAGAS} are satisfied and the theorem follows. \hfill\IEEEQED

\end{appendices}

\bibliographystyle{IEEEtran}
\bibliography{IEEEabrv,bib_file}

\begin{IEEEbiography}%
[{\includegraphics[width=1in,height=1.25in,clip,keepaspectratio]{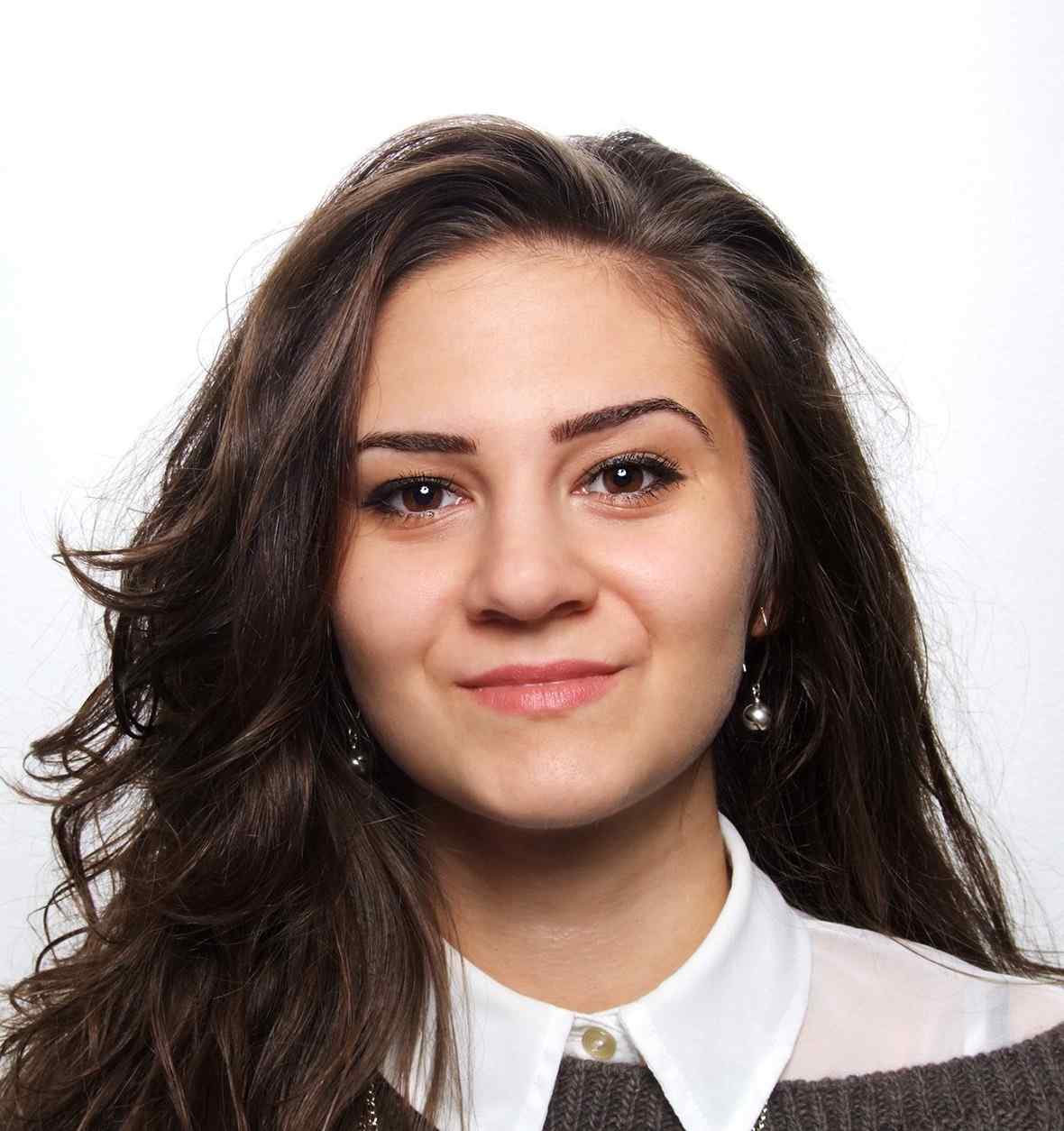}}]%
{Irina Suboti\'c} is a PhD student at the Automatic Control Laboratory at the Swiss Federal Institute of Technology (ETH) Z\"urich, Switzerland. She obtained her bachelor's degree in Electrical Engineering from the University of Belgrade, Serbia, and a master's degree in Robotics Systems and Control from ETH Z\"urich, Switzerland. She finished her thesis in collaboration with the National Renewable Energy Laboratory (NREL) in Golden, CO, USA. Her current research covers control and optimization of low-inertia power systems.
\end{IEEEbiography}

\begin{IEEEbiography}[{\includegraphics[width=1in,height=1.25in,clip,keepaspectratio]{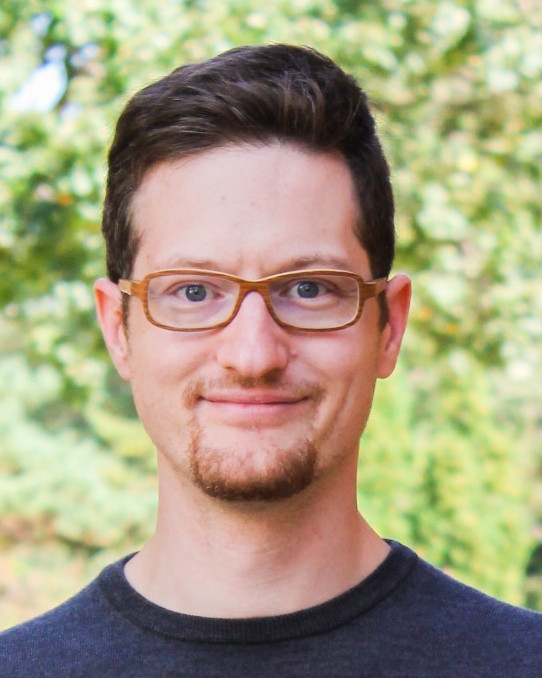}}]{Dominic Gro\ss{}} (S'11 - M'15) is an Assistant Professor with the Department of Electrical and Computer Engineering at the University of Wisconsin-Madison, Madison, WI, USA. From 2016 to 2019 he was a postdoctoral researcher at the Automatic Control Laboratory of ETH Zürich, Switzerland and he was with Volkswagen Group’s Research Division in Wolfsburg, Germany from 2014 to 2015. He received a Diploma degree in Mechatronics from the University of Kassel, Germany, in 2010, and a Ph.D. degree in Electrical Engineering from the same university in 2014. His research interests include distributed control and optimization of complex networked systems with applications in power systems dominated by power electronic devices.
\end{IEEEbiography}

\begin{IEEEbiography}%
[{\includegraphics[width=1in,height=1.25in,clip,keepaspectratio]{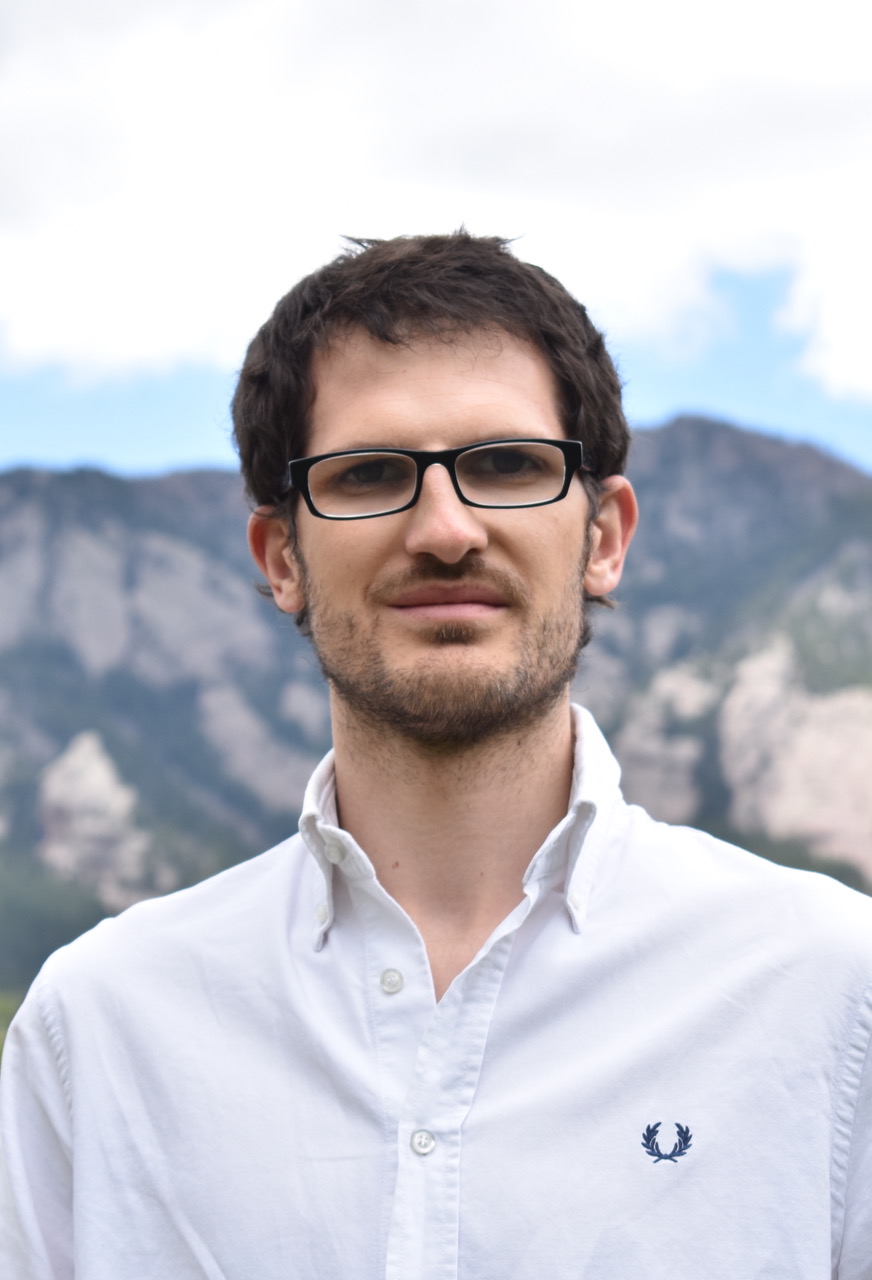}}]%
{Marcello Colombino} Marcello Colombino is a senior scientist with the Automatic Control Laboratory (IfA) at ETH, Z\"{u}rich. He received his Ph.D. degree in Electrical Engineering from IfA in 2016 and he remained as a postdoctoral researcher until 2017. in In 2018 he was with the Optimization and Control Group at the Power Systems Engineering Center of the National Renewable Energy Laboratory (NREL) in Golden, CO, USA and in 2019 he was an Assistant Professor with the Electrical and Computer Engineering Department at McGill University, Montréal, QC, Canada.  His research interests include robust and distributed control and optimization methods with application to control of energy systems, synchronization, robust control, network theory, positive systems and game-theoretical approaches to control of large-scale systems.
\end{IEEEbiography}

\begin{IEEEbiography}[{\includegraphics[width=1in,height=1.25in,clip,keepaspectratio]{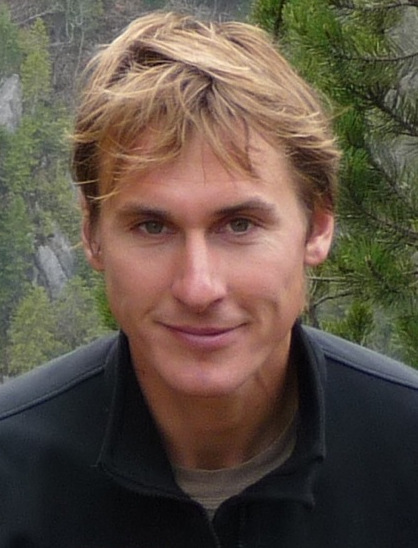}}]{Florian D\"orfler} is an Associate Professor at the Automatic Control Laboratory at ETH Zürich. He received his Ph.D. degree in Mechanical Engineering from the University of California at Santa Barbara in 2013, and a Diplom degree in Engineering Cybernetics from the University of Stuttgart in 2008. From 2013 to 2014 he was an Assistant Professor at the University of California Los Angeles. His primary research interests are centered around control, optimization, and system theory with applications in network systems such as electric power grids, robotic coordination, and social networks. He is a recipient of the distinguished young research awards by IFAC (Manfred Thoma Medal 2020) and EUCA (European Control Award 2020). His students were winners or finalists for Best Student Paper awards at the European Control Conference (2013, 2019), the American Control Conference (2016), and the PES PowerTech Conference (2017). He is furthermore a recipient of the 2010 ACC Student Best Paper Award, the 2011 O. Hugo Schuck Best Paper Award, the 2012-2014 Automatica Best Paper Award, the 2016 IEEE Circuits and Systems Guillemin-Cauer Best Paper Award, and the 2015 UCSB ME Best PhD award.
\end{IEEEbiography}

\end{document}